# Mixed variable structural optimization using mixed variable system Monte Carlo tree search formulation

**Fu-Yao Ko[1] · Katsuyuki Suzuki[1] · Kazuo Yonekura[1]**

**Abstract**

A novel method called mixed variable system Monte Carlo tree search (MVSMCTS) formulation is presented for optimization problems considering various types of variables with single and mixed continuous-discrete system. This method utilizes a reinforcement learning algorithm with improved Monte Carlo tree search (IMCTS) formulation. For sizing and shape optimization of truss structures, the design variables are the cross-sectional areas of the members and the nodal coordinates of the joints. MVSMCTS incorporates update process and accelerating technique for continuous variable and combined scheme for single and mixed system. Update process indicates that once a solution is determined by MCTS with automatic mesh generation in continuous space, it is used as the initial solution for next search tree. The search region should be expanded from the mid-point, which is the design variable for initial state. Accelerating technique is developed by decreasing the range of search region and the width of search tree based on the number of meshes during update process. Combined scheme means that various types of variables are coupled in only one search tree. Through several examples, it is demonstrated that this framework is suitable for mixed variable structural optimization. Moreover, the agent can find optimal solution in a reasonable time, stably generates an optimal design, and is applicable for practical engineering problems.

## 1 Introduction

Structural optimization has received considerable attention during the past decades due to limited material resources (Vanderplaats 1982). For sizing and shape optimization of truss structures, it becomes necessary to optimize cross-sectional areas of members and nodal coordinates of joints simultaneously (Kunar and Chan 1976; Topping 1983; Kuritz and Fleury 1989). From a weight-saving aspect, this optimization problem can provide more reduction in weight than purity sizing optimization (Kirsch and Topping 1992; Gholizadeh 2013). However, this problem is considered to be more challenging because the two types of variables involved are of fundamentally different nature.

Combining these may produce a different rate of convergence and induce the problem to be ill-conditioning (Vanderplaats and Moses 1972; Lipson and Gwin 1977; Wu and Chow 1995). To overcome the difficulties, two main classes of optimization methods have been published using mathematical programming (Yonekura and Kanno 2010; Kanno and Fujita 2018) and metaheuristics (Vargas et al. 2019; Lagaros et al. 2023).

Reinforcement learning (RL) is a branch of machine learning that aims to train an action taker called agent to take actions to maximize the cumulative numerical reward signal (Sutton and Barto 2018). Markov decision process (MDP) is the most basic theoretical model and mathematical expression for RL problems. In recent years, Q-learning-based RL approaches have been successfully used for structural design problems (Hayashi and Ohsaki 2020; Hayashi and Ohsaki 2022; Kupwiwat et al. 2024). It can be applied to structural optimization problems where it is not easy to find the optimal solutions beforehand. Moreover, RL approach does not require desired outputs as training data. However, the most significant limitation of this method is the high computational cost for training the agent. (Hayashi and Ohsaki 2020).

Recently, a Monte Carlo tree search-based RL algorithm is proposed to solve truss optimization problems. MCTS is a best-first search algorithm to solve sequential decision problems. MCTS starts with a single root node and expands the search tree that relies on random simulations rather than on full exploration. MCTS can select the optimal decision in the large search space with very little domain-specific prior knowledge. However, a huge amount of memory usage and computational resource is required to construct a search tree (Browne et al. 2012; Zuccotto et al. 2024).

A novel MCTS-based RL algorithm (Luo et al. 2022a; Luo et al. 2022b) was developed to generate optimal truss layout considering sizing, shape, and topology. Continuous design variables are considered in this study. Moreover, an RL algorithm using improved Monte Carlo tree search (IMCTS) formulation (Ko et al. 2024) was presented for discrete optimum design of truss structures. It utilizes the search tree with multiple root nodes. This formulation includes update process and accelerating technique for discrete variable, the best reward, and terminal condition. The agent is trained to minimize the weight of the truss subjected to stress and displacement constraints. The numerical results demonstrate that IMCTS formulation can find optimal solution at a low computational cost, stably generates an optimal design, and is applicable for practical engineering problems and multi-objective optimization of structures. However, IMCTS formulation is only applied to structural optimization problems considering one type of variable with discrete cases.

In this paper, mixed variable system Monte Carlo tree search (MVSMCTS) formulation based on IMCTS formulation is developed to deal with both sizing and shape variable. MVSMCTS formulation incorporates update process and accelerating technique for continuous variable and combined scheme for single and mixed continuous-discrete system. A united computational framework of update process is formed for truss design with both discrete variable and continuous one. The difference is that the solution is found by MCTS with automatic mesh generation in continuous space. Combined scheme indicates that various types of variables are considered in only one search tree. Through several numerical examples, it is demonstrated that MVSMCTS is suitable for mixed variable structural optimization. Also, the agent can find optimal solution with low computational cost, stably generates an optimal design, and is applicable for practical engineering problems.

## 2 A brief review on formulation of sizing-shape truss optimization problems and IMCTS formulation for discrete truss optimization

MVSMCTS formulation is applicable for mixed variable structural optimization. The optimal design of truss structures for sizing and shape is used as an example in this study. Therefore, this section first briefly reviews its formulation. Then, an RL algorithm using IMCTS formulation for discrete optimum design of truss structures is recalled because its concepts will be utilized in the proposed algorithm.

### 2.1 Formulation of the sizing and shape optimization of truss structures

The aim of the mixed variable structural optimization is to minimize the weight of structures considering various types of variables while satisfying design constraints. In this class of optimization problems, there are two types of design variables: sizing variable (cross-sectional areas of members) and shape variable (coordinates of nodes). The optimization problem can be described in mathematical formulations as below:

$$\boldsymbol{X} = \left(X_1, X_2, \ldots, X_{i_X}, \ldots, X_{g_X}\right) \tag{1a}$$

$$\boldsymbol{Y} = \left(Y_1, Y_2, \ldots, Y_{i_Y}, \ldots, Y_{g_Y}\right) \tag{1b}$$

Find

$$X_{i_X} \in \mathbf{D}_{X,i_X} = \left\{\left(d_{X,i_X}\right)_1, \left(d_{X,i_X}\right)_2, \ldots, \left(d_{X,i_X}\right)_h, \ldots, \left(d_{X,i_X}\right)_{b_X}\right\} \tag{1c}$$

$$Y_{i_Y} \in \mathbf{D}_{Y,i_Y} = \left\{\left(d_{Y,i_Y}\right)_1, \left(d_{Y,i_Y}\right)_2, \ldots, \left(d_{Y,i_Y}\right)_h, \ldots, \left(d_{Y,i_Y}\right)_{b_Y}\right\} \tag{1d}$$

$$X_{i_X,\min} \le X_{i_X} \le X_{i_X,\max} \tag{1e}$$

$$Y_{i_Y,\min} \le Y_{i_Y} \le Y_{i_Y,\max} \tag{1f}$$

Minimize

$$W(\boldsymbol{X}, \boldsymbol{Y}) = \rho \sum_{i_X=1}^{g_X} \left( X_{i_X} \sum_{j=1}^{m_{i_X}} L_{i_X,j}\left(Y_{i_Y}\right) \right) \tag{1g}$$

where $\boldsymbol{X}$ and $\boldsymbol{Y}$ are the sizing vector containing the cross-sectional areas and the shape vector containing the coordinates of nodes; $X_{i_X}$ is the cross-sectional area of the members belonging to group $i_X$; $Y_{i_Y}$ is the coordinate of the node $i_Y$; $X_{i_X}$ and $Y_{i_Y}$ are selected from a list of available discrete values or can be any continuous value; $g_X$ and $g_Y$ are the number of design variables for sizing and shape; $\mathbf{D}_{X,i_X}$ and $\mathbf{D}_{Y,i_Y}$ are the lists including all available discrete values arranged in ascending sequences for $X_{i_X}$ and $Y_{i_Y}$; $\left(d_{X,i_X}\right)_h$ and $\left(d_{Y,i_Y}\right)_h$ are the elements in lists $\mathbf{D}_{X,i_X}$ and $\mathbf{D}_{Y,i_Y}$; $h$ is the index of permissive discrete variable; $b_X$ and $b_Y$ are the number of available variables for $X_{i_X}$ and $Y_{i_Y}$; $X_{i_X,\min}$ and $X_{i_X,\max}$ are the lower and upper limit for $X_{i_X}$; $Y_{i_Y,\min}$ and $Y_{i_Y,\max}$ are the lower and upper limit for $Y_{i_Y}$; $W(\boldsymbol{X},\boldsymbol{Y})$ is the objective function measuring the weight of the truss; $\rho$ is the material density; $m_{i_X}$ is the number of members in group $i_X$; and $L_{i_X,j}(Y_{i_Y})$ is the length of the member $j$ in group $i_X$, which is the function of $Y_{i_Y}$. The constraints can be stated as follows:

$$\text{Subject to} \quad \begin{aligned}
&\mathbf{K}(\boldsymbol{X},\boldsymbol{Y})\boldsymbol{u} = \boldsymbol{F} && \text{(2a)}\\
&\boldsymbol{F} = (F_1^x, F_1^y, F_1^z, F_2^x, F_2^y, F_2^z, \dots, F_k^x, F_k^y, F_k^z, \dots) && \text{(2b)}\\
&\varepsilon_{i_X,j} = \boldsymbol{B}_{i_X,j}\boldsymbol{u}_{i_X,j} && \text{(2c)}\\
&\boldsymbol{B}_{i_X,j} = (-\boldsymbol{e}_{i_X,j}^{\mathrm{tr}} \quad \boldsymbol{e}_{i_X,j}^{\mathrm{tr}}) && \text{(2d)}\\
&\boldsymbol{u}_{i_X,j} = (\boldsymbol{u}_{i_X,j}^1 \quad \boldsymbol{u}_{i_X,j}^2)^{\mathrm{tr}} && \text{(2e)}\\
&\sigma_{i_X,j} = \frac{E}{L_{i_X,j}}\varepsilon_{i_X,j} && \text{(2f)}\\
&\sigma_{\min} \le \sigma_{i_X,j} \le \sigma_{\max}, i_X = 1,2,\dots,g_X, j = 1,2,\dots,m_{i_X} && \text{(2g)}\\
&\delta_{\min} \le \delta_k \le \delta_{\max}, k = 1,2,\dots,c && \text{(2h)}
\end{aligned}$$

where $\mathbf{K}(\boldsymbol{X},\boldsymbol{Y})$ is the global stiffness matrix of the entire truss; $\boldsymbol{u}$ is a vector with the displacements of the non-suppressed nodes of the entire truss; $\boldsymbol{F}$ is a vector with the given external forces at these nodes; $F_k^x$, $F_k^y$, and $F_k^z$ are the $x$, $y$, and $z$ components of external force acting on node $k$; $\varepsilon_{i_X,j}$ is the elongation of the member; $\boldsymbol{e}_{i_X,j}$ is a unit vector along the member so that it points from node 1 to node 2; $\boldsymbol{u}_{i_X,j}$ is a vector with the displacements of the end points of the member; tr in the superscript of a given vector is the transpose operator; $E$ is the modulus of elasticity; and $c$ is the number of nodes in the truss. From Eqs. 2d and 2e, the stress $\sigma_{i,j}$ in each member of the truss is compared with $\sigma_{\max}$ and $\sigma_{\min}$. The displacement $\delta_k$ of node $k$ of the truss is compared with $\delta_{\max}$ and $\delta_{\min}$. $\sigma_{\max}$ and $\sigma_{\min}$ are the maximum allowable normal stress limit for both tension and compression. $\delta_{\max}$ and $\delta_{\min}$ are the maximum and minimum permissible nodal displacement value.

Fig. 1 illustrates the example of two-bar planar truss to introduce the notation mentioned above. The number of design variables for sizing and shape are both 2 in this example. Moreover, 3 nodes are included. $x_k$, $y_k$, and $z_k$ are the coordinates of node $k$ in the $x$, $y$, and $z$ directions. $R_k^x$, $R_k^y$, and $R_k^z$ are the $x$, $y$, and $z$ components of reaction force applied to node $k$.

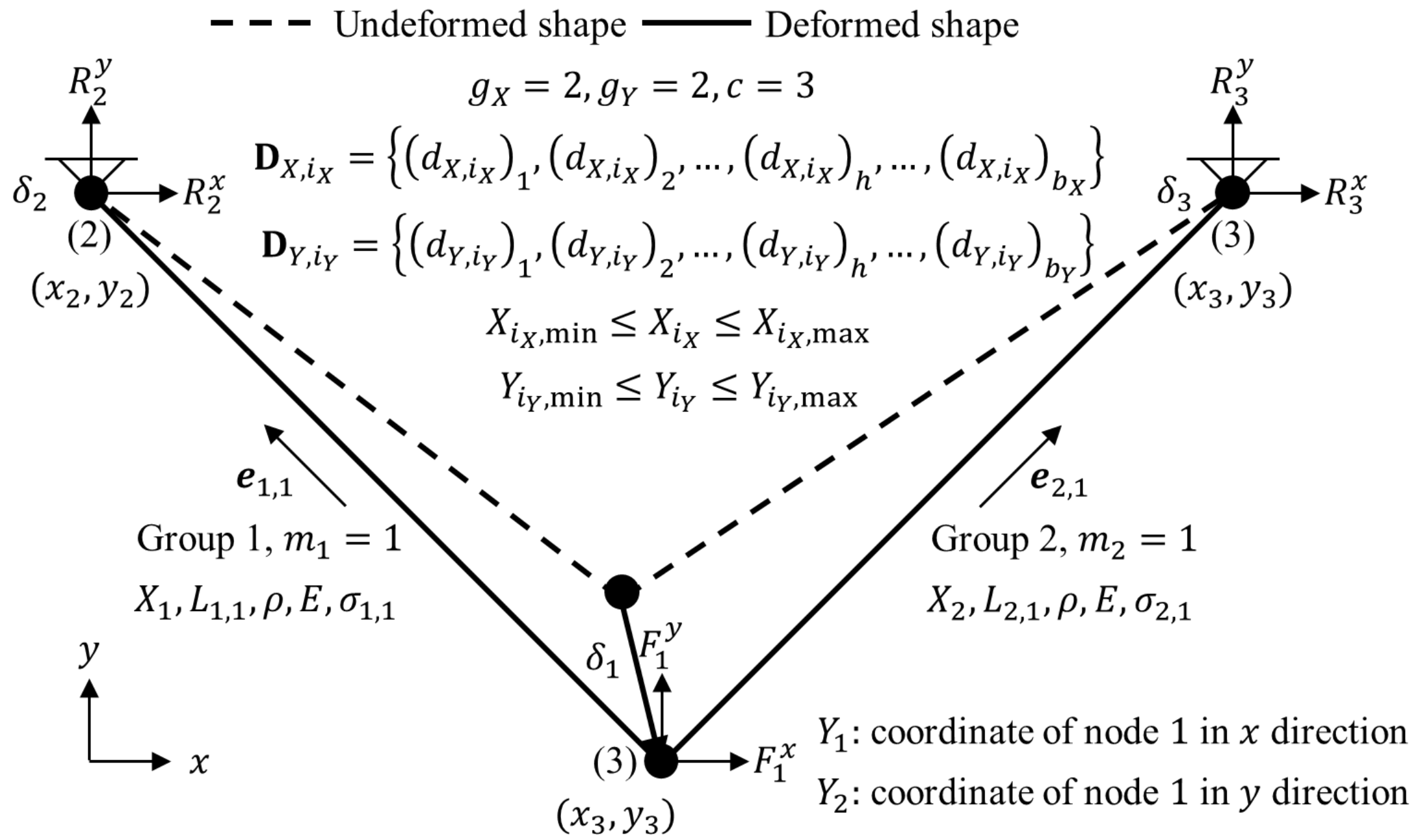

**Fig. 1** Two-bar planar truss for sizing and shape optimization of truss structures

## 2.2 IMCTS formulation for discrete sizing optimization of truss structures

An RL-based algorithm using IMCTS formulation was first proposed by Ko et al. (2024) to solve sizing optimization of truss structures with discrete variable. IMCTS formulation proposed the novel concept including update process and accelerating technique for discrete variable, the best reward, and terminal condition. These concepts are used in MVSMCTS formulation, which will be further explained in the next section.

### 2.2.1 Update process for discrete sizing variable and the best reward

The four main components of an MDP for discrete sizing optimization of truss structures were developed by Ko et al. (2024). A state is represented as a set of numerical data of nodes and members. Action is defined as determining the cross-sectional area of each member from a list of discrete values. The environment deterministically leads to a unique next state after executing action from one state. The reward is assigned to the terminal state, which is computed as follows:

$$r_T = (\alpha / W_T)^2 \tag{3}$$

where $r_T$ is the reward for terminal state $T$, which is a dimensionless quantity; $\alpha$ is the minimum weight; and $W_T$ is the weight of the truss for terminal state.

A search tree is constructed in each round for IMCTS formulation. MCTS originates from the root node denoted as the initial state of the search tree. The four strategic steps of MCTS are repeated until some predefined computational budget is reached (Browne et al. 2012). In the selection step, the upper confidence bound (UCB) for IMCTS formulation is defined as follows:

$$U_I = V_I + C\sqrt{\frac{\ln N}{n_I}} \tag{4}$$

where $I$ is the index of the node; $U_I$, $V_I$, and $n_I$ are the UCB, the estimate of state value, and the number of times the node $I$ has been visited, respectively; $N$ is the total number of simulations executed from the parent node; and $C$ is the constant parameter that adjusts the selection strategy. For IMCTS formulation, the best reward is used in the backpropagation step, which is calculated as follows:

$$V_I \leftarrow \max\left(V_I, G_{\tau_N}\right) \tag{5}$$

where $G_{\tau_N}$ is the simulation result of the path $\tau_N$ from the root node or the parent node. Moreover, $n_I$ increases by 1 for all nodes along the path.

MCTS starts from the root node or the parent node and continues executing the four steps. After reaching maximum number of iterations for that node, a child node with the largest estimate of state value is selected. This process is called policy improvement. Then, the selected node is regarded as the parent node for the next policy improvement. After many policy improvements, terminal node is determined. Then, the final state $\overline{s^p}$, the final weight $\overline{W^p}$, and the final sizing vector $\overline{\boldsymbol{X}^p}$ in round $p$ are determined.

Update process shown in Fig. 2 indicates that $\overline{\boldsymbol{X}^p}$ is used as the sizing vector $\boldsymbol{X}_0^{p+1}$ for initial state $s_0^{p+1}$ in round $(p+1)$. The starting point of the search for sizing variable is the maximum value. Therefore, all sizing variables for initial state $s_0^1$ in round 1 is equal to the largest element in a list $\mathbf{D}_{X,i_X}$, i.e., $d_{X,b_X}$. In order to describe the action and action space of initial and intermediate state for discrete sizing variable, a list $\mathbf{D}_{X,i_X}^p$ belonging to $\mathbf{D}_{X,i_X}$ in round $p$ is defined as follows:

$$\mathbf{D}_{X,i_X}^p = \left\{\left(d_{X,i_X}^p\right)_1, \dots, \left(d_{X,i_X}^p\right)_h, \dots, \left(d_{X,i_X}^p\right)_\mu, \dots, \left(d_{X,i_X}^p\right)_{\beta_X^p}\right\} \tag{6}$$

where $\left(d_{X,i_X}^p\right)_h$ is the element in a list $\mathbf{D}_{X,i_X}^p$; $\left(d_{X,i_X}^p\right)_\mu$ is the median of a list $\mathbf{D}_{X,i_X}^p$ and equal to $\left(X_{i_X}^p\right)_0$; and $\beta_X^p$ is the number of elements in a list $\mathbf{D}_{X,i_X}^p$.

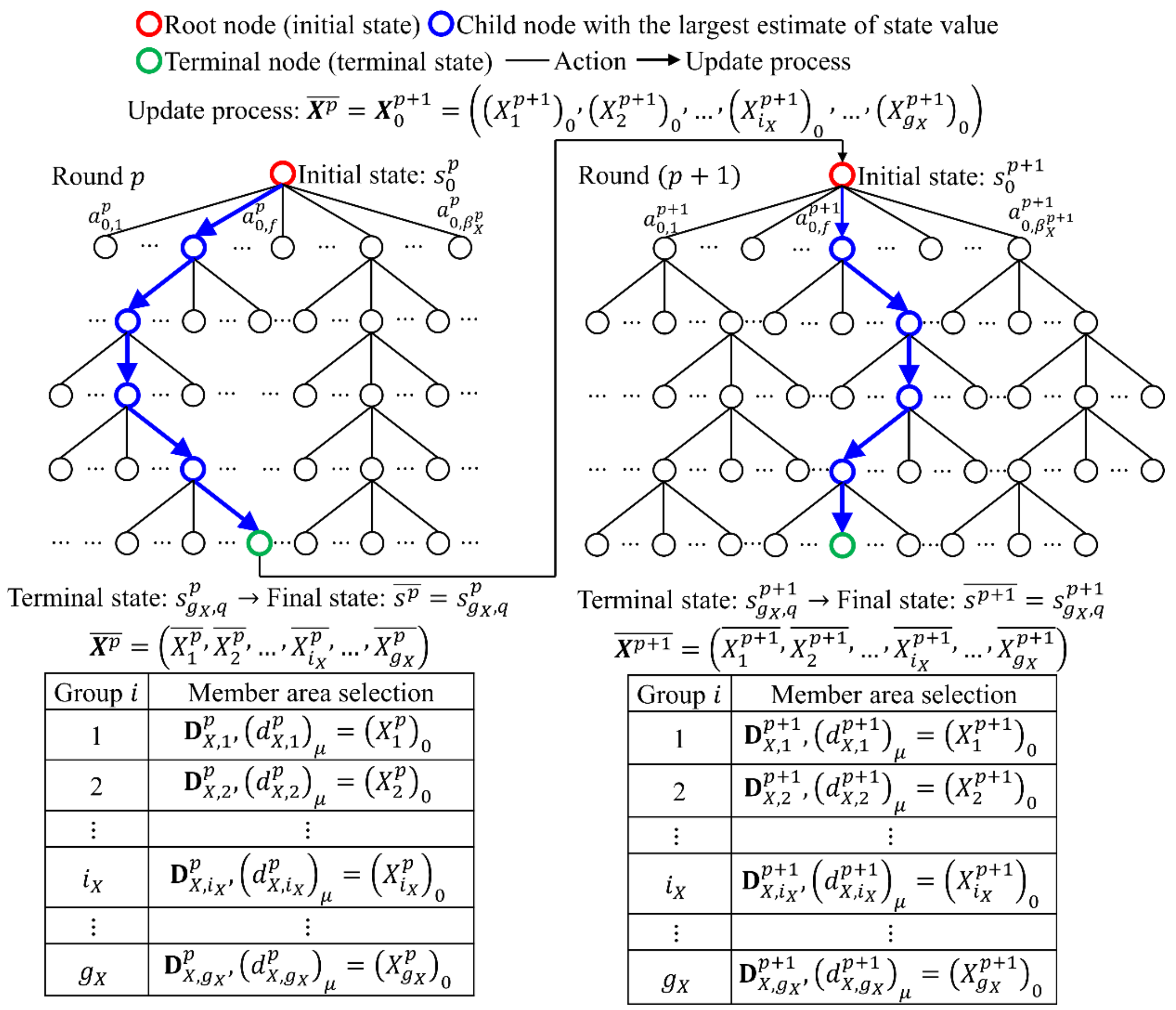

**Fig. 2** Update process for discrete sizing variable

### 2.2.2 Accelerating technique for discrete sizing variable

Accelerating technique for discrete sizing variable is developed by decreasing the width of search tree during the update process. The width of search tree is based on $\beta_X^p$. Three types of accelerating techniques are considered in the IMCTS formulation: (1) geometric decay, (2) linear decrease, and (3) step reduction. Geometric decay is the most efficient approach, which is computed as follows:

$$\beta_X^1 = \begin{cases} b_X \text{ if } b_X \text{ is odd number} \\ b_X + 1 \text{ if } b_X \text{ is even number} \end{cases} \tag{7a}$$

$$\varphi_X^p = \beta_X^1 \times \gamma_X^{\left\lceil \frac{p-1}{\epsilon_X} \right\rceil} \tag{7b}$$

$$\phi_X^p = \lfloor \varphi_X^p \rfloor \tag{7c}$$

$$\omega_X^p = \begin{cases} \phi_X^p \text{ if } \phi_X^p \text{ is odd number} \\ \phi_X^p + 1 \text{ if } \phi_X^p \text{ is even number} \end{cases} \tag{7d}$$

$$\beta_X^p = \max\left(3, \omega_X^p\right) \; (p > 1) \tag{7e}$$

where $\gamma_X$ and $\epsilon_X$ are constant parameters to adjust $\beta_X^p$, which are set to 0.5 and 3; $\left\lceil \frac{p-1}{\epsilon_X} \right\rceil$ is the least integer greater than or equal to $\frac{p-1}{\epsilon_X}$; $\lfloor \varphi_X^p \rfloor$ is the greatest integer less than or equal to $\varphi_X^p$; and $\max(3, \omega_X^p)$ is the maximum value between 3 and $\omega_X^p$.

### 2.2.3 Terminal condition

For terminal condition, improvement factor $\eta$ and counter $\theta$ are used to ensure the convergence of the algorithm. The improvement factor is defined as follows:

$$\eta = |(\overline{W^p} - \min(\mathbf{S}))/\min(\mathbf{S}) \times 100\%| \tag{8}$$

where $\mathbf{S}$ is a list to store the final weight $\overline{W^p}$; and min $(\mathbf{S})$ is the smallest element in a list $\mathbf{S}$. Before executing an algorithm, $\overline{W^0}$ is in a list $\mathbf{S}$, and $\theta$ is set to 0. $\overline{W^0}$ is the maximum weight. At the end of the round, the improvement factor is calculated, and then $\overline{W^p}$ is inserted in a list $\mathbf{S}$. When $\eta < \eta_{\min}$, $\theta$ increases by 1. When $\theta \geq \theta_{\max}$, the algorithm terminates. $\eta_{\min}$ and $\theta_{\max}$ are the critical value of improvement factor and the maximum number of counters for termination, which are set to 0.01% and 3. At this time, min $(\mathbf{S})$ is the optimal solution of the optimization problem.

## 3 MVSMCTS formulation for mixed variable structural optimization-Take sizing and shape optimization of truss structures as an example

Update process for discrete variable outlined in the previous section can be easily applied for continuous variable. Then, update process for discrete shape variable is briefly mentioned, which is similar to that for discrete sizing variable. Finally, MVSMCTS formulation considering different types of variables with mixed continuous-discrete system is described. Sizing and shape variable are considered in this study.

### 3.1 MVSMCTS formulation for continuous sizing and shape variable

### 3.1.1 Update process

For update process, it is desirable to handle the sizing and shape variable in the same computational framework. It means that $\overline{\boldsymbol{X}^p}$ and $\overline{\boldsymbol{Y}^p}$ in round $p$ is utilized as the sizing vector $\boldsymbol{X}_0^{p+1}$ and the shape vector $\boldsymbol{Y}_0^{p+1}$ for initial state $s_0^{p+1}$ in round $(p+1)$. The starting point of the search for sizing variable is the maximum value as described in Sect. 2.2.1. However, the starting point of a path for shape variable is the center of entire search region. Therefore, $\boldsymbol{X}_0^1$ and $\boldsymbol{Y}_0^1$ is expressed as follows:

$$\boldsymbol{X}_0^1 = \left((X_1^1)_0, \dots, \left(X_{i_X}^1\right)_0, \dots, (X_{g_X}^1)_0\right) = (X_{1,\max}, \dots, X_{i_X,\max}, \dots, X_{g_X,\max}) \tag{9}$$

$$\boldsymbol{Y}_0^1 = \left((Y_1^1)_0, \dots, \left(Y_{i_Y}^1\right)_0, \dots, (Y_{g_Y}^1)_0\right) = (Y_{1,\mathrm{med}}, \dots, Y_{i_Y,\mathrm{med}}, \dots, Y_{g_Y,\mathrm{med}}) \tag{10a}$$

$$Y_{i_Y,\mathrm{med}} = 0.5 \times (Y_{i_Y,\min} + Y_{i_Y,\max}) \tag{10b}$$

The search region in each round is determined from the design variable for initial state. Therefore, the left-hand side and right-hand side of search region for sizing and shape variable in round $p$ are defined as follows:

$$X_{i_X,\min}^p = \left(X_{i_X}^p\right)_0 - A_X \times \xi_{X,i_X}^p \tag{11a}$$

$$X_{i_X,\max}^p = \left(X_{i_X}^p\right)_0 + (1 - A_X) \times \xi_{X,i_X}^p \tag{11b}$$

$$Y_{i_Y,\min}^p = \left(Y_{i_Y}^p\right)_0 - A_Y \times \xi_{Y,i_Y}^p \tag{12a}$$

$$Y_{i_Y,\max}^p = \left(Y_{i_Y}^p\right)_0 + (1 - A_Y) \times \xi_{Y,i_Y}^p \tag{12b}$$

where $X_{i_X,\min}^p$ and $Y_{i_Y,\min}^p$ are the left-hand side of search region; $X_{i_X,\max}^p$ and $Y_{i_Y,\max}^p$ are the right-hand side of search region; $\xi_{X,i_X}^p$ and $\xi_{Y,i_Y}^p$ are the range of search region, which decreases when more round is reached; and $A_X$ and $A_Y$ are constant parameters to determine the position in search region for expansion. The design variable for initial state is set to the left-hand side, the center, and the right-hand side of search region when $A_X$ and $A_Y$ are all equal to 0, 0.5, and 1. The schematic diagram for expansion of search region is depicted in Fig. 3.

In order to adopt MCTS in continuous space, uniform meshes are generated automatically in search region in each round. Therefore, lists $\mathbf{Y}_{X,i_X}^p$ and $\mathbf{Y}_{Y,i_Y}^p$ for continuous sizing and shape variable in round $p$ are as follows:

$$\mathbf{Y}_{X,i_X}^p = \left\{\left(v_{X,i_X}^p\right)_1, \left(v_{X,i_X}^p\right)_2, \dots, \left(v_{X,i_X}^p\right)_h, \dots, \left(v_{X,i_X}^p\right)_{\kappa_X^p}\right\} \tag{13}$$

$$\mathbf{Y}_{Y,i_Y}^p = \left\{\left(v_{Y,i_Y}^p\right)_1, \left(v_{Y,i_Y}^p\right)_2, \dots, \left(v_{Y,i_Y}^p\right)_h, \dots, \left(v_{Y,i_Y}^p\right)_{\kappa_Y^p}\right\} \tag{14}$$

where $\left(v_{X,i_X}^p\right)_h$ and $\left(v_{Y,i_Y}^p\right)_h$ are the elements in lists $\mathbf{Y}_{X,i_X}^p$ and $\mathbf{Y}_{Y,i_Y}^p$; and $\kappa_X^p$ and $\kappa_Y^p$ are the number of elements in lists $\mathbf{Y}_{X,i_X}^p$ and $\mathbf{Y}_{Y,i_Y}^p$. It is worth noting that $\left(v_{X,i_X}^p\right)_1$, $\left(v_{X,i_X}^p\right)_{\kappa_X^p}$, $\left(v_{Y,i_Y}^p\right)_1$, and $\left(v_{Y,i_Y}^p\right)_{\kappa_Y^p}$ are $X_{i_X,\min}^p$, $X_{i_X,\max}^p$, $Y_{i_Y,\min}^p$, and $Y_{i_Y,\max}^p$, respectively. For continuous sizing and shape variable, action $a_{l_X,f}^p$ $(l_X = 0, 1, \dots, g_X - 1)$ and $a_{l_Y,f}^p$ $(l_Y = 0, 1, \dots, g_Y - 1)$ for initial and intermediate state is expressed as follows:

$$a_{l_X,f}^p\text{: } M_{X,i_X}^p\text{: } 1 \rightarrow 0 \text{ and } X_{i_X}^p\text{: } \left(X_{i_X}^p\right)_0 \rightarrow \left(v_{X,i_X}^p\right)_h \ (i_X = l_X + 1, f = h) \tag{15}$$

$$a^p_{l_Y,f}: M^p_{Y,i_Y}: 1 \to 0 \text{ and } Y^p_{i_Y}: \left(Y^p_{i_Y}\right)_0 \to \left(v^p_{Y,i_Y}\right)_h \ (i_Y = l_Y + 1, f = h) \tag{16}$$

where $M^p_{X,i_X}$ and $M^p_{Y,i_Y}$ are used to confirm when the cross-sectional area of group $i_X$ and the nodal coordinate of node $i_Y$ are determined in round $p$; $X^p_{i_X}$ is the cross-sectional area; and $Y^p_{i_Y}$ is the nodal coordinate. $M^p_{X,i_X} = 1$ and $M^p_{Y,i_Y} = 1$ when $X_{i_X}$ and $Y_{i_Y}$ are not determined, and $M^p_{X,i_X} = 0$ and $M^p_{Y,i_Y} = 0$ otherwise. The action spaces of $s^p_0$, $s^p_{l_X,q}$, and $s^p_{l_Y,q}$ are as follows:

$$\boldsymbol{\mathcal{A}}\left(s^p_0\right) = \left\{a^p_{0,1}, a^p_{0,2}, \dots, a^p_{0,f}, \dots, a^p_{0,\kappa^p_X}\right\} \tag{17a}$$

$$\boldsymbol{\mathcal{A}}\left(s^p_{l_X,q}\right) = \left\{a^p_{l_X,1}, a^p_{l_X,2}, \dots, a^p_{l_X,f}, \dots, a^p_{l_X,\kappa^p_X}\right\} \tag{17b}$$

$$\boldsymbol{\mathcal{A}}\left(s^p_0\right) = \left\{a^p_{0,1}, a^p_{0,2}, \dots, a^p_{0,f}, \dots, a^p_{0,\kappa^p_Y}\right\} \tag{18a}$$

$$\boldsymbol{\mathcal{A}}\left(s^p_{l_Y,q}\right) = \left\{a^p_{l_Y,1}, a^p_{l_Y,2}, \dots, a^p_{l_Y,f}, \dots, a^p_{l_Y,\kappa^p_Y}\right\} \tag{18b}$$

The schematic diagram of update process for continuous sizing and shape variable is depicted in Fig. 3.

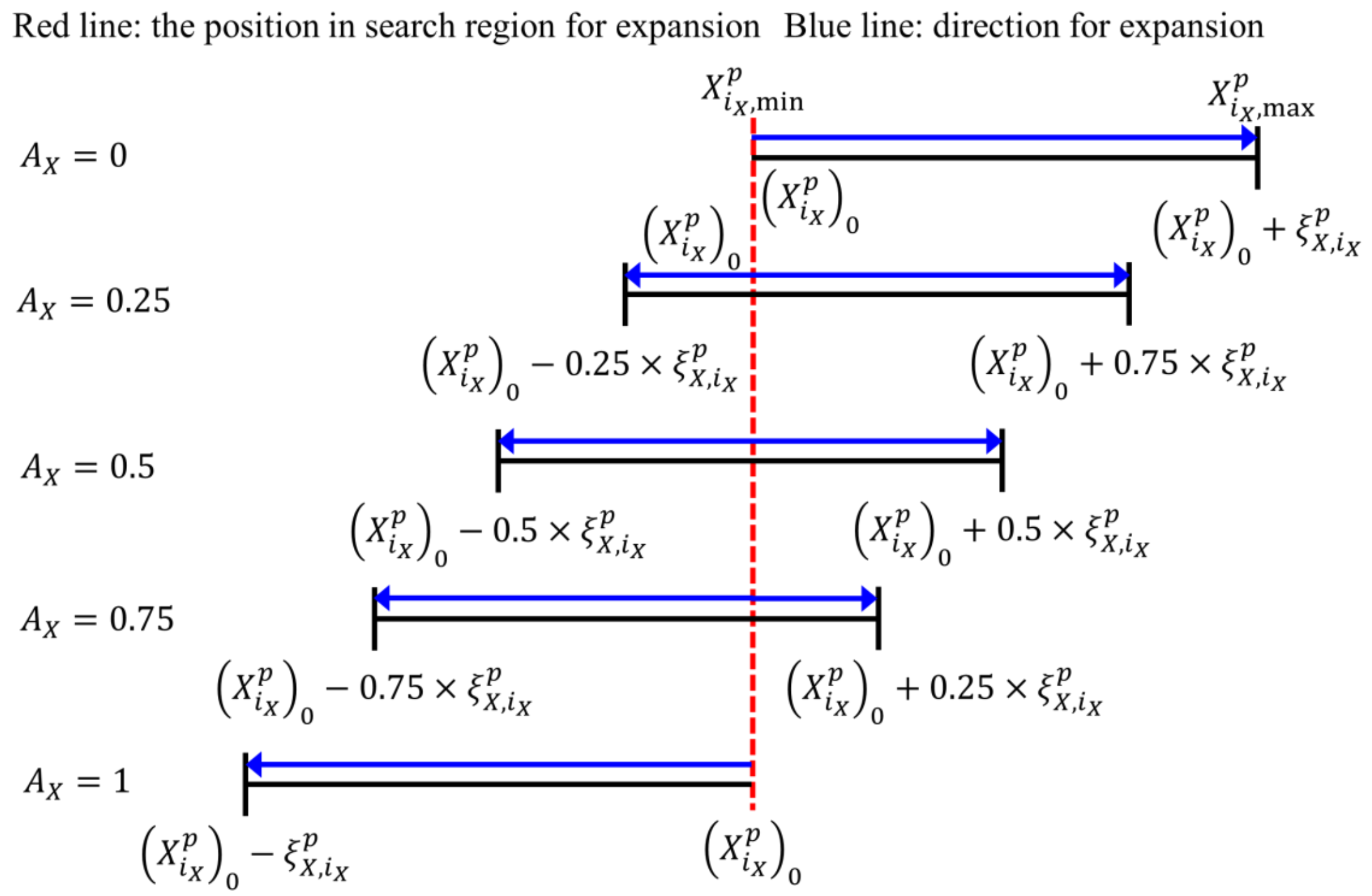

**Fig. 3** Schematic diagram for expansion of search region in continuous space

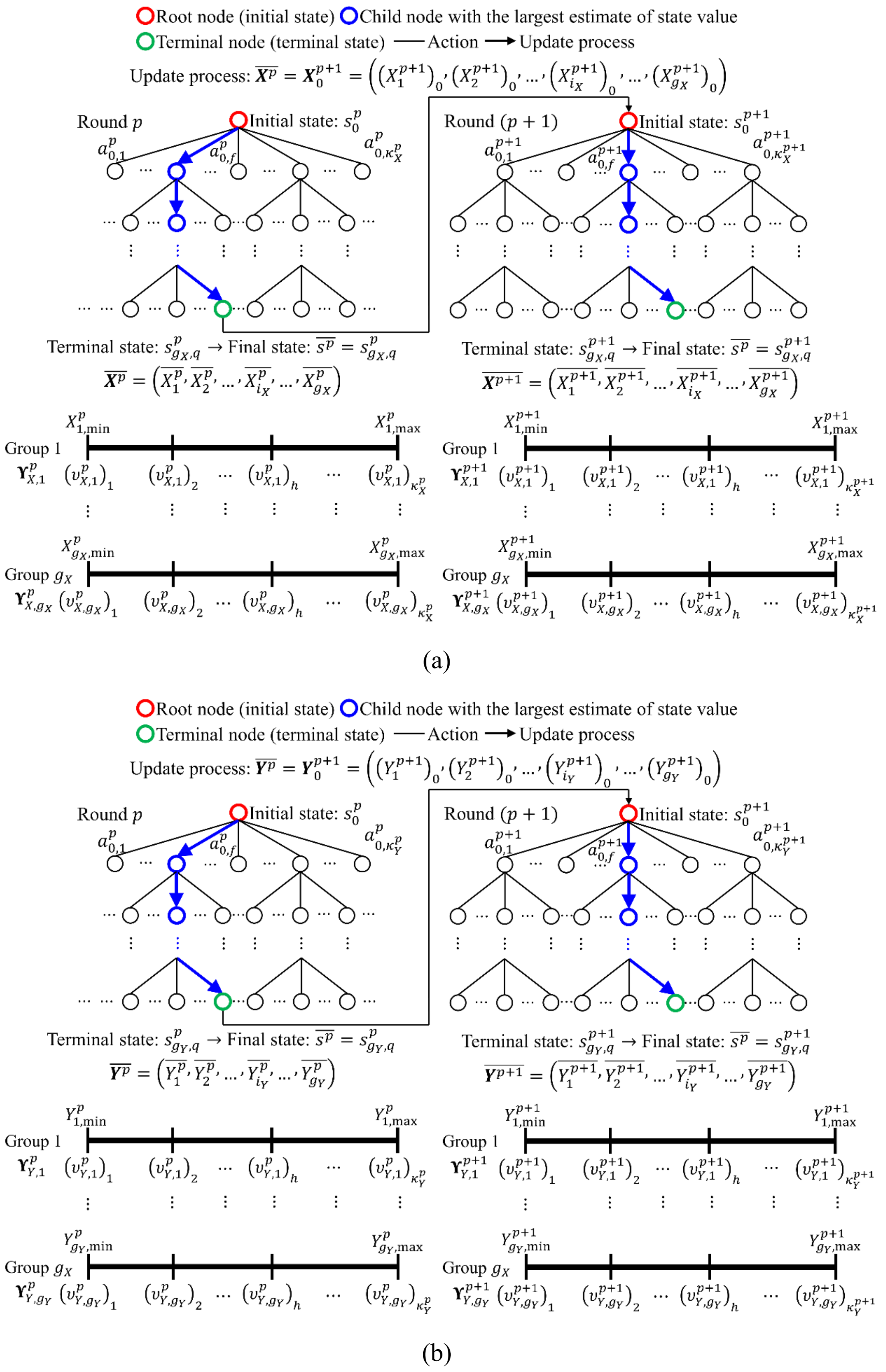

**Fig. 4** Schematic diagram of update process for continuous (a) sizing and (b) shape variable

### 3.1.2 Accelerating technique

The range of search region and the width of search tree are not changed during update process when accelerating technique is not considered, as shown in Fig. 5(a). It means that $\xi^p_{X,i_X}$, $\kappa^p_X$, $\xi^p_{Y,i_Y}$, and $\kappa^p_Y$ are constants in each round, which is formulated as follows:

$$\xi^1_{X,i_X} = \xi^2_{X,i_X} = \cdots = \xi^p_{X,i_X} = \cdots = X_{i_X,\max} - X_{i_X,\min} \tag{19a}$$
$$\kappa^1_X = \kappa^2_X = \cdots = \kappa^p_X = \cdots = g_X \tag{19b}$$

$$\xi^1_{Y,i_Y} = \xi^2_{Y,i_Y} = \cdots = \xi^p_{Y,i_Y} = \cdots = Y_{i_Y,\max} - Y_{i_Y,\min} \tag{20a}$$
$$\kappa^1_Y = \kappa^2_Y = \cdots = \kappa^p_Y = \cdots = g_Y \tag{20b}$$

From Eqs. 19 and 20, it is shown that $\xi^p_{X,i_X}$ and $\xi^p_{Y,i_Y}$ are equal to the difference between the lower and the upper limit for sizing and shape. Moreover, $\kappa^p_X$ and $\kappa^p_Y$ are equal to the number of design variables for sizing and shape.

Accelerating technique is developed by decreasing the range of search region and the width of search tree as the update process proceeds, as shown in Fig. 5(b). Three types of accelerating techniques mentioned in Sect. 2.2.2 are used to determine the range of search region. For geometric decay, $\xi^p_{X,i_X}$ and $\xi^p_{Y,i_Y}$ are computed as follows:

$$\xi^1_{X,i_X} = X_{i_X,\max} - X_{i_X,\min} \tag{21a}$$
$$\xi^p_{X,i_X} = \xi^1_{X,i_X} \times \lambda_{X,geo}^{\left\lceil \frac{p-1}{\varpi_{X,geo}} \right\rceil} \tag{21b}$$

$$\xi^1_{Y,i_Y} = Y_{i_Y,\max} - Y_{i_Y,\min} \tag{22a}$$
$$\xi^p_{Y,i_Y} = \xi^1_{Y,i_Y} \times \lambda_{Y,geo}^{\left\lceil \frac{p-1}{\varpi_{Y,geo}} \right\rceil} \tag{22b}$$

where $\lambda_{X,geo}$, $\varpi_{X,geo}$, $\lambda_{Y,geo}$, and $\varpi_{Y,geo}$ are constant parameters for geometric decay to adjust $\xi^p_{X,i_X}$ and $\xi^p_{Y,i_Y}$ in each round; and $\left\lceil \frac{p-1}{\varpi_{X,geo}} \right\rceil$ and $\left\lceil \frac{p-1}{\varpi_{Y,geo}} \right\rceil$ are the least integers greater than or equal to $\frac{p-1}{\varpi_{X,geo}}$ and $\frac{p-1}{\varpi_{Y,geo}}$. In this study, $\lambda_{X,geo}$ and $\lambda_{Y,geo}$ are all set to 0.5, and $\varpi_{X,geo}$ and $\varpi_{Y,geo}$ are set to 3 and 5 in this study. As seen from Eqs. 21a and 22a, $\xi^1_{X,i_X}$ and $\xi^1_{Y,i_Y}$ are equal to the difference between the minimum and the maximum allowable value for sizing and shape. In order to reduce the execution time, $\xi^p_{X,i_X}$ and $\xi^p_{Y,i_Y}$ are decreasing geometrically when reaching more rounds. Hence, Eqs. 21b and 22b are utilized to satisfy this requirement.

For linear decrease, $\xi^p_{X,i_X}$ and $\xi^p_{Y,i_Y}$ are calculated as follows:

$$\xi^1_{X,i_X} = X_{i_X,\max} - X_{i_X,\min} \tag{23a}$$
$$\xi^p_{X,i_X} = \xi^1_{X,i_X} - \lambda_{X,lin}(p-1) \tag{23b}$$
$$\xi^p_{X,i_X} = \max\left(\xi^{lin}_{X,i_X}, \xi^p_{X,i_X}\right)\ (p > 1) \tag{23c}$$

$$\xi_{X,i_X}^{lin} = \lambda_{X,lin} \times \xi_{X,i_X}^{1} \tag{23d}$$

$$\xi_{Y,i_Y}^{1} = Y_{i_Y,\max} - Y_{i_Y,\min} \tag{24a}$$
$$\xi_{Y,i_Y}^{p} = \xi_{Y,i_Y}^{1} - \lambda_{Y,lin}(p-1) \tag{24b}$$
$$\xi_{Y,i_Y}^{p} = \max\left(\xi_{Y,i_Y}^{lin}, \xi_{Y,i_Y}^{p}\right)\ (p > 1) \tag{24c}$$
$$\xi_{Y,i_Y}^{lin} = \lambda_{Y,lin} \times \xi_{Y,i_Y}^{1} \tag{24d}$$

where $\lambda_{X,lin}$ and $\lambda_{Y,lin}$ are constant parameters for linear decrease to adjust $\xi_{X,i_X}^{p}$ and $\xi_{Y,i_Y}^{p}$ in each round, which are all set to 0.05 in this study; and $\xi_{X,i_X}^{lin}$ and $\xi_{Y,i_Y}^{lin}$ are the critical values of linear decrease for $\xi_{X,i_X}^{p}$ and $\xi_{Y,i_Y}^{p}$, respectively. From Eqs. 23b and 24b, it is seen that $\xi_{X,i_X}^{p}$ and $\xi_{Y,i_Y}^{p}$ are deceasing linearly when more round is reached. Eqs. 23c and 24c are employed to ensure that $\xi_{X,i_X}^{p}$ and $\xi_{Y,i_Y}^{p}$ are all positive real numbers during update process.

For step reduction, $\xi_{X,i_X}^{p}$ and $\xi_{Y,i_Y}^{p}$ can be expressed as follows:

$$\xi_{X,i_X}^{1} = X_{i_X,\max} - X_{i_X,\min} \tag{25a}$$
$$\xi_{X,i_X}^{p} = \xi_{X,i_X}^{1} - \lambda_{X,red}\left\lfloor\frac{p-1}{\varpi_{X,red}}\right\rfloor \tag{25b}$$
$$\xi_{X,i_X}^{p} = \max\left(\xi_{X,i_X}^{red}, \xi_{X,i_X}^{p}\right)\ (p > 1) \tag{25c}$$
$$\xi_{X,i_X}^{red} = \lambda_{X,red} \times \xi_{X,i_X}^{1} \tag{25d}$$

$$\xi_{Y,i_Y}^{1} = Y_{i_Y,\max} - Y_{i_Y,\min} \tag{26a}$$
$$\xi_{Y,i_Y}^{p} = \xi_{Y,i_Y}^{1} - \lambda_{Y,red}\left\lfloor\frac{p-1}{\varpi_{Y,red}}\right\rfloor \tag{26b}$$
$$\xi_{Y,i_Y}^{p} = \max\left(\xi_{Y,i_Y}^{red}, \xi_{Y,i_Y}^{p}\right)\ (p > 1) \tag{26c}$$
$$\xi_{Y,i_Y}^{red} = \lambda_{Y,red} \times \xi_{Y,i_Y}^{1} \tag{26d}$$

where $\lambda_{X,red}$, $\varpi_{X,red}$, $\lambda_{Y,red}$, $\varpi_{Y,red}$ are constant parameters for step reduction to adjust $\xi_{X,i_X}^{p}$ and $\xi_{Y,i_Y}^{p}$ in each round, which are set to 0.05, 3, 0.05, and 5 in this study; and $\xi_{X,i_X}^{red}$ and $\xi_{Y,i_Y}^{red}$ are the critical values of step reduction for $\xi_{X,i_X}^{p}$ and $\xi_{Y,i_Y}^{p}$, respectively. In Eqs. 25b and 26b, it is shown that $\xi_{X,i_X}^{p}$ and $\xi_{Y,i_Y}^{p}$ are the step function of $p$.

The width of search tree is based on the number of meshes automatically generated in continuous space. In other words, $\kappa_X^p$ and $\kappa_Y^p$ are utilized to determine the width of search tree in each round. There are three types of accelerating techniques to calculate $\kappa_X^p$ and $\kappa_Y^p$, which is similar to that for discrete variable presented in Sect. 2.2.2. For geometric decay, $\kappa_X^p$ and $\kappa_Y^p$ are expressed as follows:

$$\kappa_X^1 = g_X \tag{27a}$$
$$\psi_X^p = \kappa_X^1 \times \varrho_{X,geo}^{\left\lfloor\frac{p-1}{\varsigma_{X,geo}}\right\rfloor} \tag{27b}$$
$$\kappa_X^p = \max\left(3, \left\lfloor\psi_X^p\right\rfloor\right)\ (p > 1) \tag{27c}$$

$$\kappa_Y^1 = g_Y \tag{28a}$$
$$\psi_Y^p = \kappa_Y^1 \times \varrho_{Y,geo}^{\left\lceil \frac{p-1}{\varsigma_{Y,geo}} \right\rceil} \tag{28b}$$
$$\kappa_Y^p = \max\left(3, \lfloor \psi_Y^p \rfloor\right)\ (p > 1) \tag{28c}$$

where $\varrho_{X,geo}$, $\varsigma_{X,geo}$, $\varrho_{Y,geo}$, and $\varsigma_{Y,geo}$ are constant parameters to adjust $\kappa_X^p$ and $\kappa_Y^p$ in each round, which are set to 0.5, 3, 0.5, and 5 in this study; $\left\lceil \frac{p-1}{\varsigma_{X,geo}} \right\rceil$ and $\left\lceil \frac{p-1}{\varsigma_{Y,geo}} \right\rceil$ are the least integers greater than or equal to $\frac{p-1}{\varsigma_{X,geo}}$ and $\frac{p-1}{\varsigma_{Y,geo}}$; $\lfloor \psi_X^p \rfloor$ and $\lfloor \psi_Y^p \rfloor$ are the greatest integers less than or equal to $\psi_X^p$ and $\psi_Y^p$. Eqs. 27a and 28a indicate that $\kappa_X^1$ and $\kappa_Y^1$ are equal to the number of design variables for sizing and shape. In order to reduce the computation time, $\kappa_X^p$ and $\kappa_Y^p$ are decreasing geometrically when more round is reached. Therefore, Eqs. 27b and 28b are employed to fulfill the requirements. It is worth mentioning that $\kappa_X^p$ and $\kappa_Y^p$ are unnecessary to be odd number for continuous variable since $\kappa_X^p$ and $\kappa_Y^p$ are utilized to determine the number of meshes automatically generated in search region. In Eqs. 27c and 28c, the minimum value of $\kappa_X^p$ and $\kappa_Y^p$ needs to be 3 because there are three parts in search region: the left-hand side, the center, and the right-hand side.

For linear decrease, $\kappa_X^p$ and $\kappa_Y^p$ are given as follows:

$$\kappa_X^1 = g_X \tag{29a}$$
$$\Gamma_X^p = \kappa_X^1 - \varrho_{X,lin}(p-1) \tag{29b}$$
$$\kappa_X^p = \max\left(3, \lfloor \Gamma_X^p \rfloor\right)\ (p > 1) \tag{29c}$$

$$\kappa_Y^1 = g_Y \tag{30a}$$
$$\Gamma_Y^p = \kappa_Y^1 - \varrho_{Y,lin}(p-1) \tag{30b}$$
$$\kappa_Y^p = \max\left(3, \lfloor \Gamma_Y^p \rfloor\right)\ (p > 1) \tag{30c}$$

where $\varrho_{X,lin}$ and $\varrho_{Y,lin}$ are constant parameters for linear decrease to adjust $\kappa_X^p$ and $\kappa_Y^p$ in each round, which are all set to 2 in this study. From Eqs. 29b and 30b, $\kappa_X^p$ and $\kappa_Y^p$ are decreasing linearly when reaching more rounds.

For step reduction, $\kappa_X^p$ and $\kappa_Y^p$ are expressed as follows:

$$\kappa_X^1 = g_X \tag{31a}$$
$$\Omega_X^p = \kappa_X^1 - \varrho_{X,red} \left\lfloor \frac{p-1}{\varsigma_{X,red}} \right\rfloor \tag{31b}$$
$$\kappa_X^p = \max\left(3, \lfloor \Omega_X^p \rfloor\right)\ (p > 1) \tag{31c}$$

$$\kappa_Y^1 = g_Y \tag{32a}$$
$$\Omega_Y^p = \kappa_Y^1 - \varrho_{Y,red} \left\lfloor \frac{p-1}{\varsigma_{Y,red}} \right\rfloor \tag{32b}$$
$$\kappa_Y^p = \max\left(3, \lfloor \Omega_Y^p \rfloor\right)\ (p > 1) \tag{32c}$$

where $\varrho_{X,red}$, $\varsigma_{X,red}$, $\varrho_{Y,red}$, and $\varsigma_{Y,red}$ are constant parameters for step reduction to adjust $\kappa_X^p$ and $\kappa_Y^p$ in each round, which are set to 2, 3, 2, and 5 in this study. As seen from Eqs. 31b and 32b, $\kappa_X^p$ and $\kappa_Y^p$ are the step function of $p$.

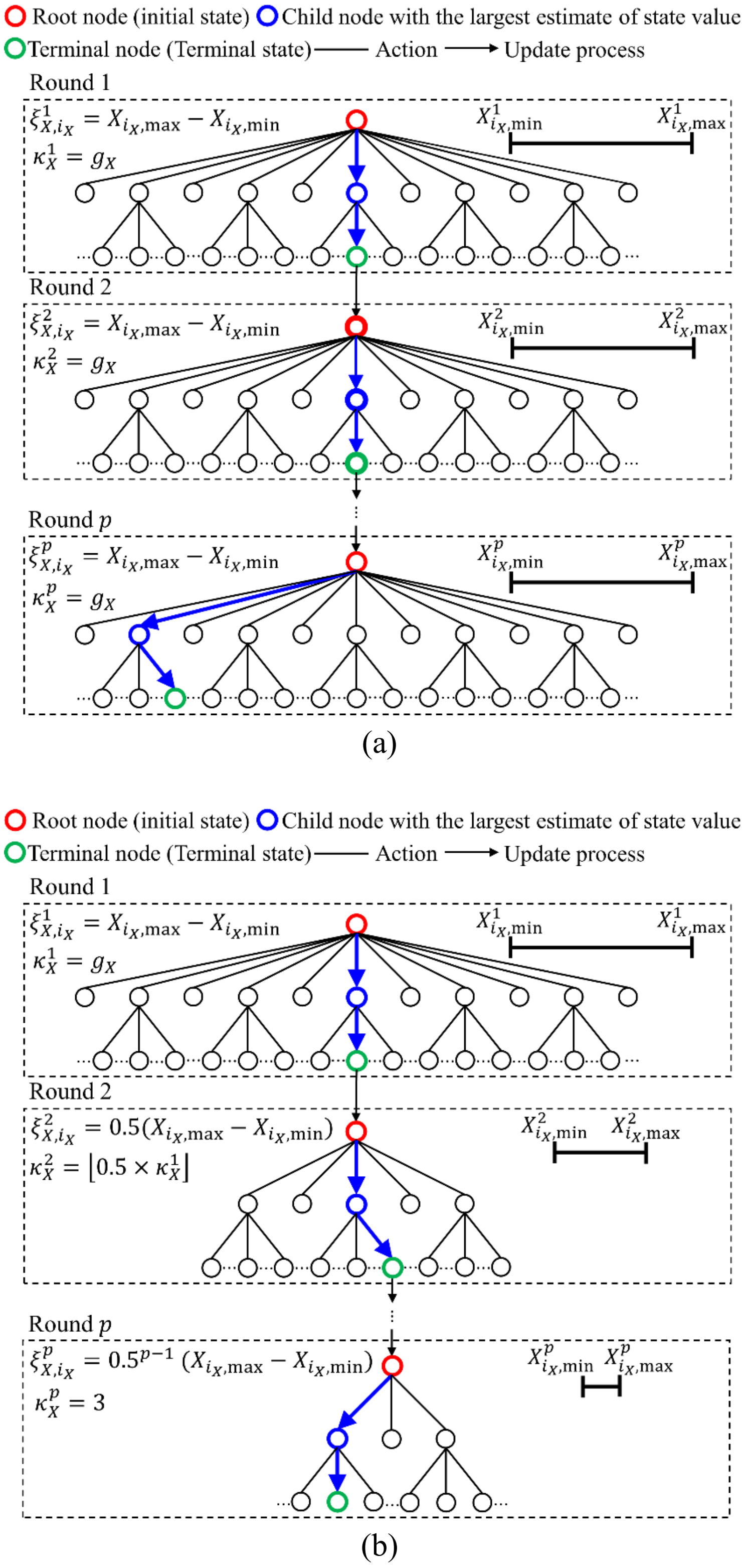

**Fig. 5** MVSMCTS formulation for continuous variable (a) without and (b) with accelerating technique (geometric decay)

### 3.2 MVSMCTS formulation for discrete shape variable

#### 3.2.1 Update process

Update process for discrete shape variable is the same as that for continuous one in Sect. 3.1.1. As mentioned previously, the starting point of the search for shape variable is the center of entire search region. Therefore, $\boldsymbol{Y}_0^1$ is described as follows:

$$\boldsymbol{Y}_0^1 = \left((Y_1^1)_0, \dots, \left(Y_{i_Y}^1\right)_0, \dots, (Y_{g_Y}^1)_0\right) = (\left(d_{Y,1}\right)_\mu, \dots, \left(d_{Y,i_Y}\right)_\mu, \dots, \left(d_{Y,g_Y}\right)_\mu) \tag{33}$$

where $\left(d_{Y,i_Y}\right)_\mu$ is the median of a list $\mathbf{D}_{Y,i_Y}$. In order to describe the action and action space of initial and intermediate state for discrete shape variable, a list $\mathbf{D}_{Y,i_Y}^p$ belonging to $\mathbf{D}_{Y,i_Y}$ in round $p$ is defined as follows:

$$\mathbf{D}_{Y,i_Y}^p = \left\{\left(d_{Y,i_Y}^p\right)_1, \dots, \left(d_{Y,i_Y}^p\right)_h, \dots, \left(d_{Y,i_Y}^p\right)_\mu, \dots, \left(d_{Y,i_Y}^p\right)_{\beta_Y^p}\right\} \tag{34}$$

where $\left(d_{Y,i_Y}^p\right)_h$ is the element in a list $\mathbf{D}_{Y,i_Y}^p$; $\left(d_{Y,i_Y}^p\right)_\mu$ is the median of a list $\mathbf{D}_{Y,i_Y}^p$ and equal to $\left(Y_{i_Y}^p\right)_0$; and $\beta_Y^p$ is the number of elements in a list $\mathbf{D}_{Y,i_Y}^p$.

#### 3.2.2 Accelerating technique

The concept of the accelerating technique for discrete shape variable is the same as that for discrete sizing variable. The width of search tree is based on $\beta_Y^p$. Three types of accelerating techniques are considered. For geometric decay, $\beta_Y^p$ is computed as follows:

$$\beta_Y^1 = \begin{cases} b_Y \text{ if } b_Y \text{ is odd number} \\ b_Y + 1 \text{ if } b_Y \text{ is even number} \end{cases} \tag{35a}$$

$$\varphi_Y^p = \beta_Y^1 \times \gamma_Y^{\left\lceil \frac{p-1}{\epsilon_Y} \right\rceil} \tag{35b}$$

$$\phi_Y^p = \lfloor \varphi_Y^p \rfloor \tag{35c}$$

$$\omega_Y^p = \begin{cases} \phi_Y^p \text{ if } \phi_Y^p \text{ is odd number} \\ \phi_Y^p + 1 \text{ if } \phi_Y^p \text{ is even number} \end{cases} \tag{35d}$$

$$\beta_Y^p = \max(3, \omega_Y^p)\ (p > 1) \tag{35e}$$

where $\gamma_Y$ and $\epsilon_Y$ are constant parameters to adjust $\beta_Y^p$, which are set to 0.5 and 3; $\left\lceil \frac{p-1}{\epsilon_Y} \right\rceil$ is the least integer greater than or equal to $\frac{p-1}{\epsilon_Y}$; $\lfloor \varphi_Y^p \rfloor$ is the greatest integer less than or equal to $\varphi_Y^p$; and $\max(3, \omega_Y^p)$ is the maximum value between 3 and $\omega_Y^p$. The formulas of linear decrease and step reduction is the same for discrete sizing variable. The only difference is that subscript $X$ is changed to $Y$.

## 3.3 MVSMCTS formulation for sizing and shape optimization of truss structures

### 3.3.1 MCTS methodology

For MVSMCTS formulation, the search tree is built in each round to consider sizing and shape variable at the same time. The node denotes the state of the current truss structure. $s_0^p$, $s_{l_{MVS},q}^p$ $(l_{MVS} = 1, 2, \dots, g_X + g_Y - 1)$, and $s_{g_X+g_Y,q}^p$ represent the initial state, the intermediate state, and the final state in round $p$, respectively. The mathematical expressions of $s_0^p$, $s_{l_{MVS},q}^p$, and $s_{g_X+g_Y,q}^p$ are as follows:

$$\text{State } s_0^p \quad \begin{array}{c} M_{X,1}^p = M_{X,2}^p = \cdots = M_{X,g_X}^p = 1 \\ M_{Y,1}^p = M_{Y,2}^p = \cdots = M_{Y,g_Y}^p = 1 \end{array} \tag{36a}$$

$$\begin{array}{r} \text{State } s_{l_{MVS},q}^p \\ (l_{MVS} = 1, 2, \dots, g_Y - 1) \end{array} \quad \begin{array}{c} M_{X,1}^p = M_{X,2}^p \dots = M_{X,g_X}^p = 1 \\ M_{Y,1}^p = M_{Y,2}^p = \cdots = M_{Y,i_Y}^p = 0 \\ M_{Y,i_Y+1}^p = M_{Y,i_Y+2}^p = \cdots = M_{Y,g_Y}^p = 1 \ (i_Y = l_{MVS}) \end{array} \tag{36b}$$

$$\begin{array}{r} \text{State } s_{l_{MVS},q}^p \\ (l_{MVS} = g_Y) \end{array} \quad \begin{array}{c} M_{X,1}^p = M_{X,2}^p = \cdots = M_{X,g_X}^p = 1 \\ M_{Y,1}^p = M_{Y,2}^p = \cdots = M_{Y,g_Y}^p = 0 \end{array} \tag{36c}$$

$$\begin{array}{r} \text{State } s_{l_{MVS},q}^p \\ (l_{MVS} = g_Y + 1, \dots, \\ g_X + g_Y - 1) \end{array} \quad \begin{array}{c} M_{X,1}^p = M_{X,2}^p = \cdots = M_{X,i_X}^p = 0 \\ M_{X,i_X+1}^p = M_{X,i_X+2}^p = \cdots = M_{X,g_X}^p = 1 \\ M_{Y,1}^p = M_{Y,2}^p = \cdots = M_{Y,g_Y}^p = 0 \ (g_Y + i_X = l_{MVS}) \end{array} \tag{36d}$$

$$\text{State } s_{g_X+g_Y,q}^p \quad \begin{array}{c} M_{X,1}^p = M_{X,2}^p = \cdots = M_{X,g_X}^p = 0 \\ M_{Y,1}^p = M_{Y,2}^p = \cdots = M_{Y,g_Y}^p = 0 \end{array} \tag{36e}$$

where $l_{MS}$ is the layer number for MVSMCTS formulation. The root node represents the initial state $s_0^p$. Nodes with terminal state $s_{g_X+g_Y,q}^p$ are called terminal nodes. The tree edges indicate the possible actions. Action $a_{l_{MVS},f}^p$ $(l_{MVS} = 0, 1, 2, \dots, g_X + g_Y - 1)$ for initial and intermediate state is described as follows:

$$a_{l_{MVS},f}^p \colon M_{Y,i_Y}^p \colon 1 \to 0, Y_{i_Y}^p \colon \left(Y_{i_Y}^p\right)_0 \to \left(v_{Y,i_Y}^p\right)_h \ (i_Y = l_{MVS} + 1, f = h) \tag{37a}$$

$$a_{l_{MVS},f}^p \colon M_{Y,i_Y}^p \colon 1 \to 0, Y_{i_Y}^p \colon \left(Y_{i_Y}^p\right)_0 \to \left(d_{Y,i_Y}^p\right)_h \ (i_Y = l_{MVS} + 1, f = h) \tag{37b}$$

$$a_{l_{MVS},f}^p \colon M_{X,i_X}^p \colon 1 \to 0, X_{i_X}^p \colon \left(X_{i_X}^p\right)_0 \to \left(v_{X,i_X}^p\right)_h \ (i_X = l_{MVS} + 1 - g_Y, f = h) \tag{37c}$$

$$a_{l_{MVS},f}^p \colon M_{X,i_X}^p \colon 1 \to 0, X_{i_X}^p \colon \left(X_{i_X}^p\right)_0 \to \left(d_{X,i_X}^p\right)_h \ (i_X = l_{MVS} + 1 - g_Y, f = h) \tag{37d}$$

where $X_{i_X}^p$ and $Y_{i_Y}^p$ are the cross-sectional area of the group $i_X$ and the nodal coordinate of the node $i_Y$ for that state. The action spaces of $s_0^p$ and $s_{l_{MVS},q}^p$ are as follows:

$$\boldsymbol{\mathcal{A}}\left(s_0^p\right) = \left\{a_{0,1}^p, \dots, a_{0,f}^p, \dots, a_{0,\kappa_Y^p}^p\right\} \tag{38a}$$

$$\boldsymbol{\mathcal{A}}\left(s_0^p\right) = \left\{a_{0,1}^p, \dots, a_{0,f}^p, \dots, a_{0,\beta_Y^p}^p\right\} \tag{38b}$$

$$\mathcal{A}(s^p_{l_{MVS},q}) = \left\{a^p_{l_{MVS},1}, \dots, a^p_{l_{MVS},f}, \dots, a^p_{l_{MVS},\kappa^p_Y}\right\} \ (l_{MVS} = 1, \dots, g_Y - 1) \tag{38c}$$

$$\mathcal{A}(s^p_{l_{MVS},q}) = \left\{a^p_{l_{MVS},1}, \dots, a^p_{l_{MVS},f}, \dots, a^p_{l_{MVS},\beta^p_Y}\right\} \ (l_{MVS} = 1, \dots, g_Y - 1) \tag{38d}$$

$$\mathcal{A}(s^p_{l_{MVS}}) = \left\{a^p_{l_{MVS},1}, \dots, a^p_{l_{MVS},f}, \dots, a^p_{l_{MVS},\kappa^p_X}\right\} \ (l_{MVS} = g_Y, \dots, g_Y + g_X - 1) \tag{38e}$$

$$\mathcal{A}(s^p_{l_{MVS},q}) = \left\{a^p_{l_{MVS},1}, \dots, a^p_{l_{MVS},f}, \dots, a^p_{l_{MVS},\beta^p_X}\right\} \ (l_{MVS} = g_Y, \dots, g_Y + g_X - 1) \tag{38f}$$

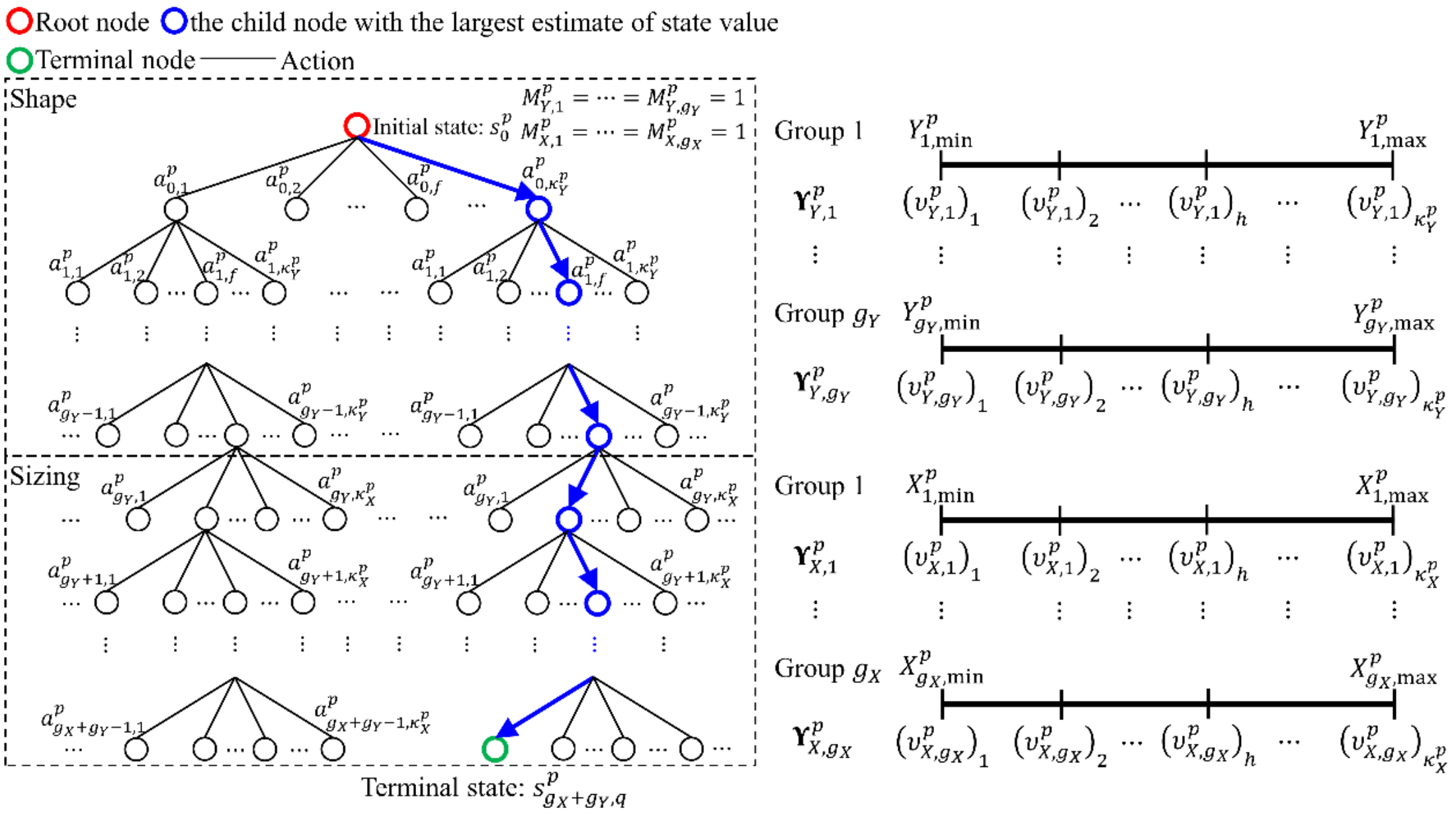

(a)

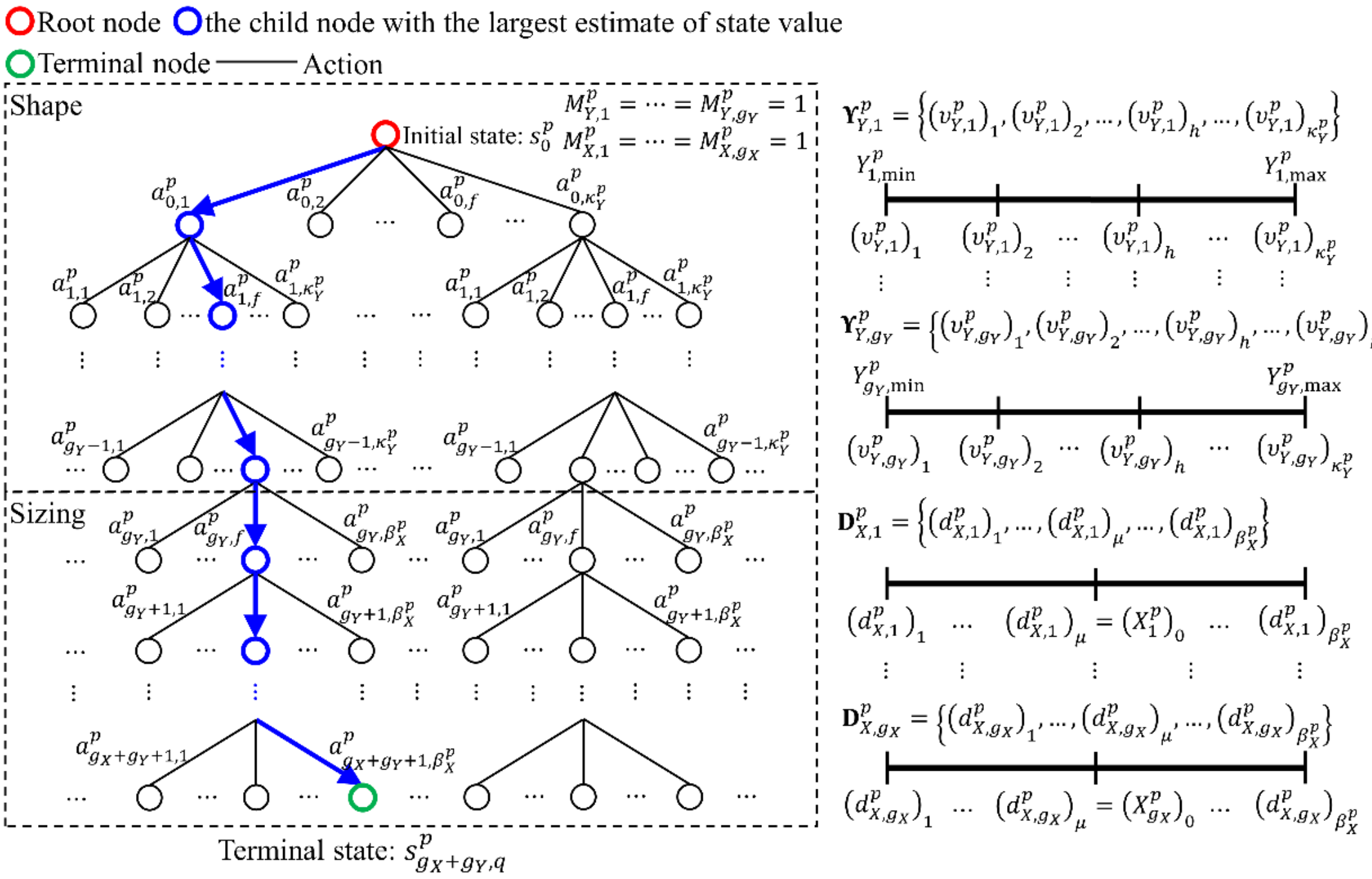

(b)

**Fig. 6** Search tree used in the MVSMCTS formulation for (a) single continuous-continuous system and (b) mixed continuous-discrete system

MCTS starts constructing the search tree with the root node following the four main steps. The UCB in Eq. 4 and the best reward in Eq. 5 for IMCTS formulation are used in the selection and the backpropagation step for MVMCTS formulation. The search tree considering sizing and shape variable with single and mixed system is shown in Fig. 6.

### 3.3.2 Policy improvement

Policy improvement for truss optimization on sizing and shape is identical to that for discrete sizing optimization outlined in Sect. 2.2.1. Sizing and shape variable are optimized in one search tree in each round. After many policy improvement steps, terminal node is reached. Then, the final state $\overline{s^p}$, the final weight $\overline{W^p}$, the final sizing vector $\overline{\boldsymbol{X}^p}$, and the final shape vector $\overline{\boldsymbol{Y}^p}$ are all determined. The schematic diagram for policy improvement is shown in Fig. 7.

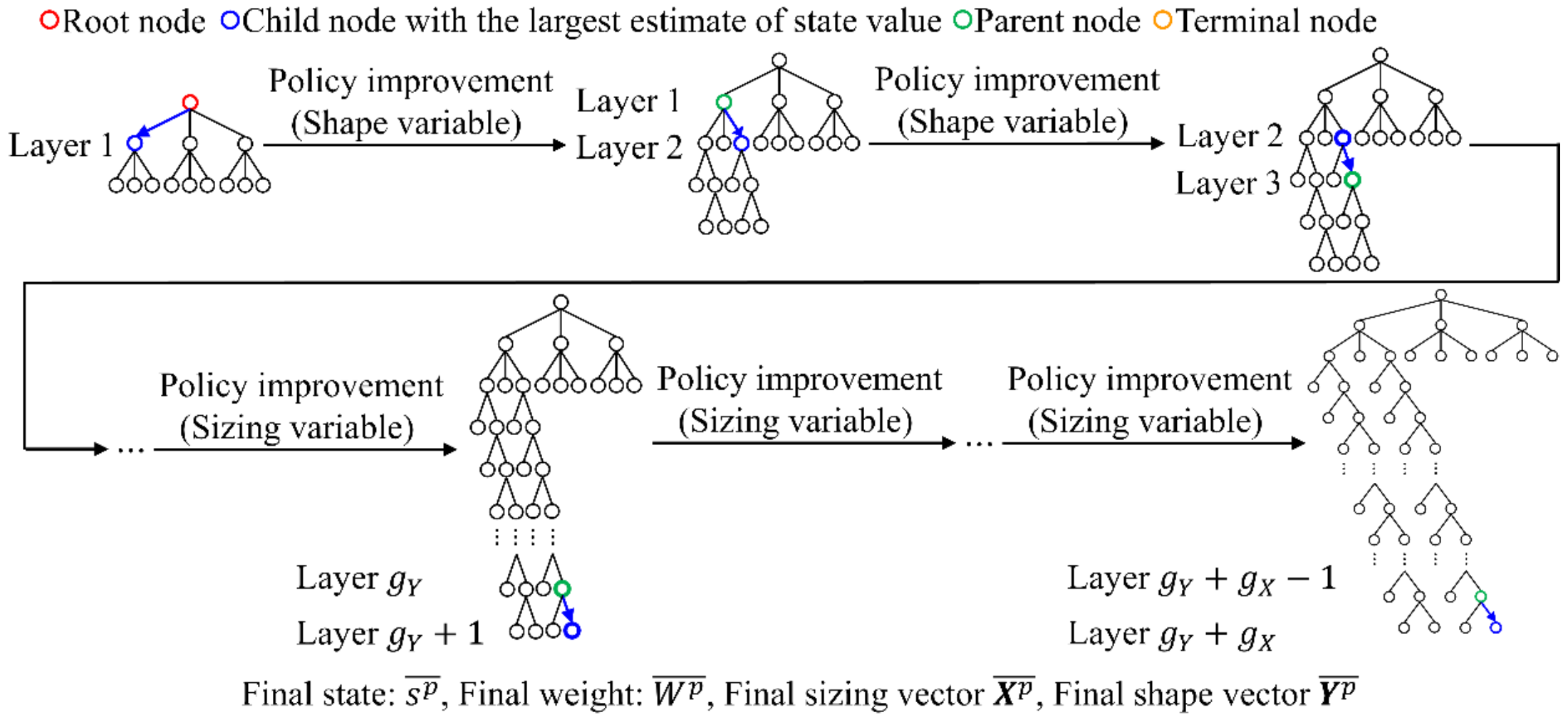

**Fig. 7** Policy improvement for MVSMCTS formulation

### 3.3.3 Update process

Update process for MVSMCTS formulation shown in Fig. 8 indicates that $\overline{\boldsymbol{X}^p}$ and $\overline{\boldsymbol{Y}^p}$ determined in round $p$ are used as the sizing vector $\boldsymbol{X}_0^{p+1}$ and shape vector $\boldsymbol{Y}_0^{p+1}$ for initial state $s_0^{p+1}$ in round $(p+1)$. All discrete and continuous sizing variables for initial state in round 1 are equal to $\left(d_{X,i_X}\right)_{b_X}$ and $X_{i_X,\max}$ shown in Sect. 2.2.1 and Eq. 9, respectively. Moreover, all discrete and continuous shape variables for initial state in round 1 are equal to $\left(d_{Y,i_Y}\right)_{b_Y}$ and $Y_{i_Y,\mathrm{med}}$ as in Eqs. 33 and 10, respectively.

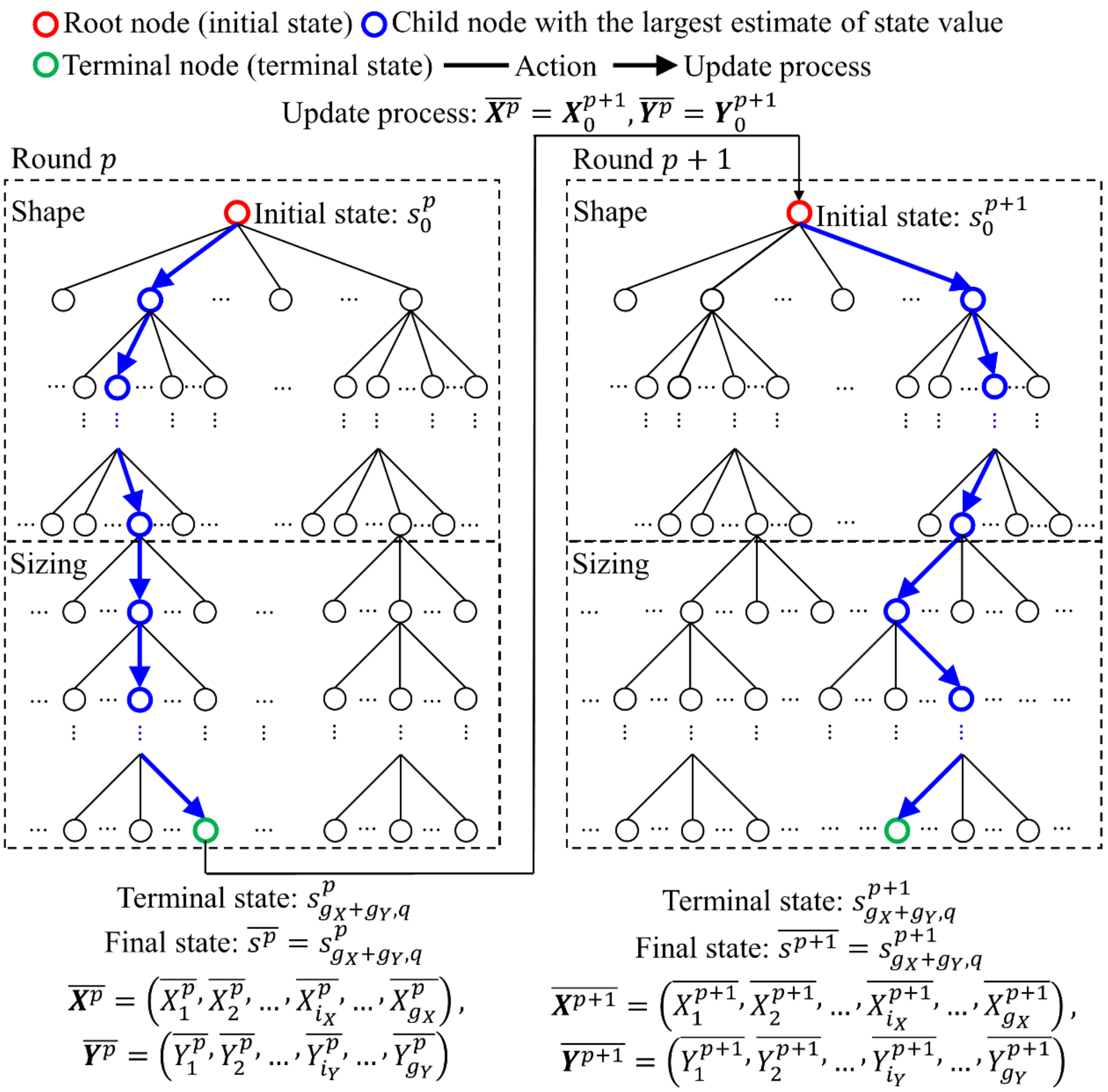

**Fig. 8** Update process for MVSMCTS formulation

Accelerating technique for MVSMCTS formulation is based on Eq. 7 for discrete sizing variable, Eqs. 21 and 27 for continuous sizing variable, Eq. 35 for discrete shape variable, and Eqs. 22 and 28 for continuous shape variable. For MVSMCTS formulation, the maximum number of iterations for root node $J_{0,\max}^p$ and parent node $J_{l_{MVS},\max}^p$ for layer $l_{MVS}$ ($l_{MVS} = 1, 2, \dots, g_X + g_Y - 1$) is expressed as follows:

For root node State $s_0^p$

$$J_{0,\max}^p = J_{MVS} \times \lceil \log_{10}(\Psi_X{}^{g_X} \times \Psi_Y{}^{g_Y}) \rceil \quad (39a)$$

For parent node State $s_{l_{MVS},q}^p$ ($l_{MVS} = 1, \dots, g_Y - 1$)

$$J_{l_{MVS},\max}^p = J_{MVS} \times \left\lceil \log_{10}\left(\Psi_X{}^{g_X} \times \Psi_Y{}^{g_Y - l_{MVS}}\right) \right\rceil \quad (39b)$$

For parent node State $s_{l_{MVS},q}^p$ ($l_{MVS} = g_Y$)

$$J_{l_{MVS},\max}^p = J_{MVS} \times \lceil \log_{10}(\Psi_X{}^{g_X}) \rceil \quad (39c)$$

For parent node State $s_{l_{MVS},q}^p$

$$J_{l_{MVS},\max}^p = J_{MVS} \times \left\lceil \log_{10}\left(\Psi_X{}^{g_X - (l_{MVS} - g_Y)}\right) \right\rceil \quad (39d)$$

$$(l_{MVS} = g_Y + 1, \ldots,$$
$$g_X + g_Y - 1)$$

where $J_{MVS}$ is a constant parameter for $J_{0,\max}^{p}$ and $J_{l_{MVS},\max}^{p}$ and is set to 3 times the number of design variables in the optimization problem in Sect. 2.1, i.e., $3(g_X + g_Y)$ in this study; $\Psi_X$ are equal to $\beta_X^p$ and $\kappa_X^p$ for discrete and continuous sizing variable; and $\Psi_Y$ are equal to $\beta_Y^p$ and $\kappa_Y^p$ for discrete and continuous shape variable. The terminal condition for IMCTS formulation outlined in Sect. 2.2.3 is also applicable for MVSMCTS formulation. The pseudo-code and the flowchart of the MVSMCTS formulation for mixed variable structural optimization are illustrated in Figs. 9 and 10.

Algorithm 1: MVSMCTS formulation for sizing and shape optimization of truss structures

**Initialize** Material density $\rho$, the modulus of elasticity $E$, the lists $\mathbf{D}_{X,i_X}$ and $\mathbf{D}_{Y,i_Y}$, the lower and upper limit $X_{i_X,\min}$, $Y_{i_Y,\min}$, $X_{i_X,\max}$, and $Y_{i_Y,\max}$<br>
**Initialize** Round number $p = 1$<br>
**Initialize** Counter $\theta = 0$ and the maximum number of counters for termination $\theta_{\max} = 3$<br>
**Initialize** Critical value of improvement factor $\eta_{\min} = 0.01$<br>
**Initialize** List $\mathbf{S} = \{\overline{W^0}\}$ and the maximum weight $\overline{W^0}$<br>
**Initialize** Design variable vector for initial state $s_0^1$ in round 1 for sizing and shape: $\boldsymbol{X}_0^1$ and $\boldsymbol{Y}_0^1$<br>
**While** $\theta \le \theta_{\max}$ **do**<br>
  Start from a search tree with only root node (initial state $s_0^p$)<br>
  Define a list $\mathbf{D}_{X,i_X}^p = \left\{\left(d_{X,i_X}^p\right)_1, \ldots, \left(d_{X,i_X}^p\right)_h, \ldots, \left(d_{X,i_X}^p\right)_\mu, \ldots, \left(d_{X,i_X}^p\right)_{\beta_X^p}\right\}$ belonging to $\mathbf{D}_{X,i_X}$ for discrete sizing variable<br>
  Define a list $\mathbf{D}_{Y,i_Y}^p = \left\{\left(d_{Y,i_Y}^p\right)_1, \ldots, \left(d_{Y,i_Y}^p\right)_h, \ldots, \left(d_{Y,i_Y}^p\right)_\mu, \ldots, \left(d_{Y,i_Y}^p\right)_{\beta_Y^p}\right\}$ belonging to $\mathbf{D}_{Y,i_Y}$ for discrete shape variable<br>
  Define a list $\mathbf{Y}_{X,i_X}^p = \left\{\left(v_{X,i_X}^p\right)_1, \left(v_{X,i_X}^p\right)_2, \ldots, \left(v_{X,i_X}^p\right)_h, \ldots, \left(v_{X,i_X}^p\right)_{\kappa_X^p}\right\}$ for continuous sizing variable, where $\left(v_{X,i_X}^p\right)_1 = X_{i_X,\min}^p$ and $\left(v_{X,i_X}^p\right)_{\kappa_X^p} = X_{i_X,\max}^p$<br>
  Define a list $\mathbf{Y}_{Y,i_Y}^p = \left\{\left(v_{Y,i_Y}^p\right)_1, \left(v_{Y,i_Y}^p\right)_2, \ldots, \left(v_{Y,i_Y}^p\right)_h, \ldots, \left(v_{Y,i_Y}^p\right)_{\kappa_Y^p}\right\}$ for continuous shape variable, where $\left(v_{Y,i_Y}^p\right)_1 = Y_{i_Y,\min}^p$ and $\left(v_{Y,i_Y}^p\right)_{\kappa_Y^p} = Y_{i_Y,\max}^p$<br>
  **While** Terminal node is reached **do**<br>
    **if** Not all shape variable are optimized **then**<br>
      **While** Maximum number of iterations is reached **do**<br>
        Repeat the four strategic steps of MCTS from the root node or the parent node<br>
        Policy improvement (shape variable): select a child node with the largest estimate of state value<br>
        This child node is regarded as the parent node<br>
    **if** Not all sizing variable are optimized **then**<br>
      **While** Maximum number of iterations is reached **do**<br>
        Repeat the four strategic steps of MCTS from the parent node<br>
        Policy improvement (sizing variable): select a child node with the largest estimate of state value<br>
        This child node is regarded as the parent node<br>
  The final state $\overline{s^p}$ is determined, which is the state of the terminal node<br>
  Sizing and shape vector $\overline{\boldsymbol{X}^p}$ and $\overline{\boldsymbol{Y}^p}$ are all determined<br>
  The final weight $\overline{W^p}$ is determined<br>
  Calculate improvement factor $\eta = |(\overline{W^p} - \min(\mathbf{S}))/\min(\mathbf{S}) \times 100\%|$<br>
  **if** $\eta \le \eta_{\min}$ **then**<br>
    $\theta \leftarrow \theta + 1$<br>
  $\overline{W^p}$ is inserted in a list $\mathbf{S}$<br>
  Update process: $\boldsymbol{X}_0^{p+1} = \overline{\boldsymbol{X}^p}$ and $\boldsymbol{Y}_0^{p+1} = \overline{\boldsymbol{Y}^p}$<br>
  $p \leftarrow p + 1$<br>
$\min(\mathbf{S})$ is the optimal solution of the sizing-shape truss optimization problems in Sect. 2.1

**Fig. 9** Pseudo-code for the MVSMCTS formulation

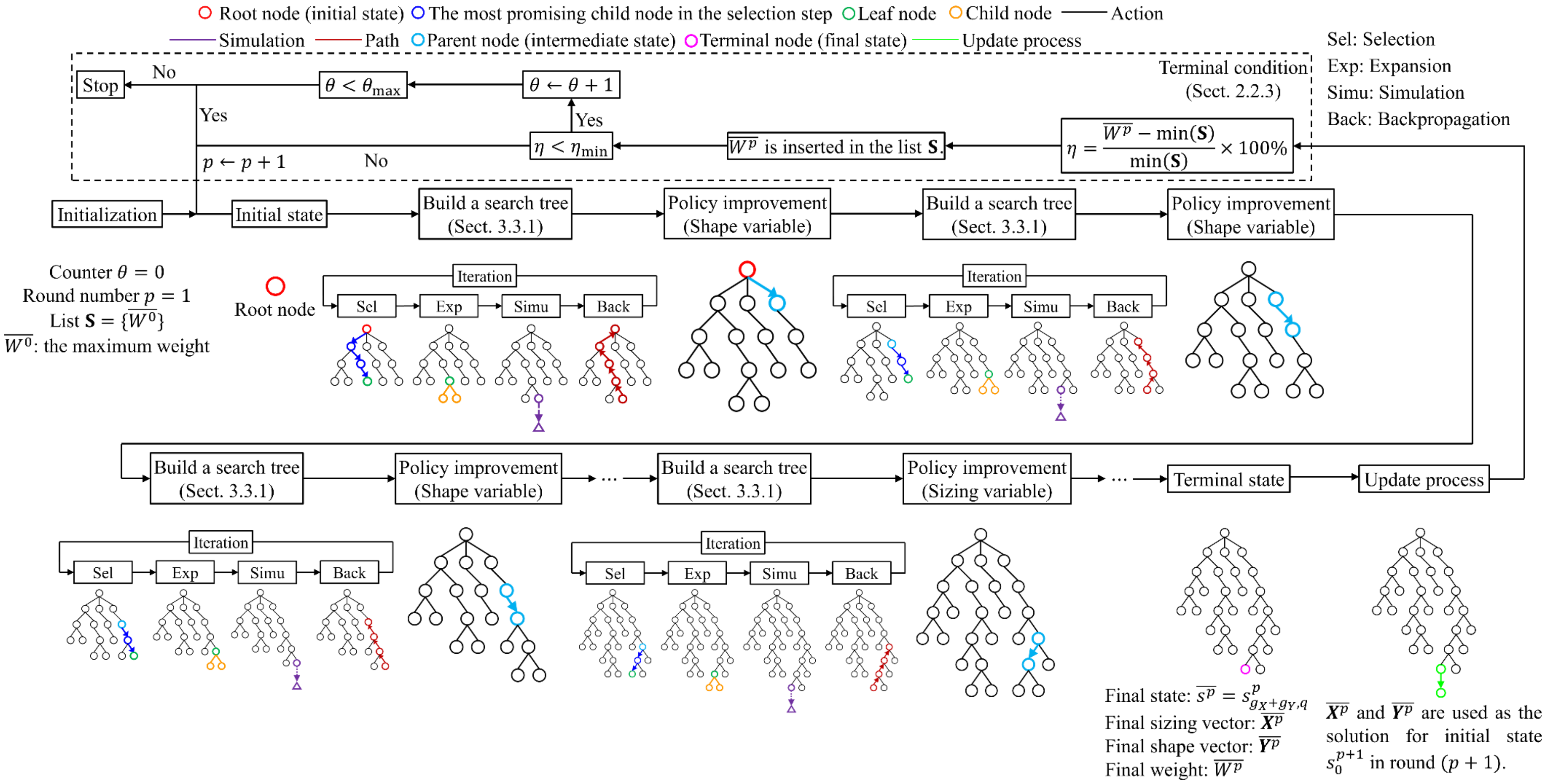

**Fig. 10** Flowchart of the MVSMCTS formulation

### 3.4 Comparisons of IMCTS and MVSMCTS formulation

Table 1 lists the characteristics of IMCTS and MVSMCTS formulation. IMCTS and MVSMCTS formulation are applicable to single type of variable with discrete system and variable types of variables with single and mixed system. Moreover, it is seen that MDP framework, four steps of MCTS, UCB, the best reward, policy improvement, update process and accelerating technique for discrete variable and terminal condition proposed in IMCTS formulation can also be used in MVSMCTS formulation.

**Table 1** Comparisons of the IMCTS and MVSMCTS formulation

| | IMCTS formulation (Ko et al. 2024) | MVSMCTS formulation |
|---|---|---|
| Type of optimization problem | Single type of variable with discrete system | Various types of variables with single and mixed continuous-discrete system |
| Example | Discrete sizing optimization of truss structure | Truss optimization on sizing and shape |
| Markov decision process | | |
| State | a set of numerical data of nodes and members | a set of numerical data of nodes and members |
| Action | determine discrete sizing variable | 1. determine discrete and continuous sizing variable<br>2. determine discrete and continuous shape variable |
| State transition | deterministic | deterministic |
| Reward function | $r_T = (\alpha/W_T)^2$ (terminal state) | $r_T = (\alpha/W_T)^2$ (terminal state) |
| Monte Carlo tree search | | |
| How to build a search tree | the four strategic steps of MCTS | the four strategic steps of MCTS |
| Selection step | UCB: $U_I = V_I + C\sqrt{\frac{\ln N}{n_I}}$ | UCB: $U_I = V_I + C\sqrt{\frac{\ln N}{n_I}}$ |
| Backpropagation step | the best reward: $V_I \leftarrow \max(V_I, G_{\tau_N})$ | the best reward: $V_I \leftarrow \max(V_I, G_{\tau_N})$ |
| Policy improvement | Start from root node and conduct policy improvement. After many policy improvement steps, terminal node is reached. | Start from root node and conduct policy improvement. After many policy improvement steps, terminal node is reached. |
| Algorithm | | |
| Update process | for discrete variable | for discrete and continuous variable<br>■ continuous: generate uniform meshes automatically in search region |
| Search tree | multiple root nodes | ■ multiple root nodes<br>■ consider various types of variables in only one search tree |
| Accelerating technique | for discrete variable<br>● decrease the width of search tree<br>● reduce maximum number of iterations | ■ for discrete variable<br>● decrease the width of search tree<br>● reduce maximum number of iterations<br>■ for continuous variable<br>● decrease the width of search tree based on the number of meshes<br>● decrease the range of search region<br>● reduce maximum number of iterations |
| Terminal condition | Improvement factor $\eta$ and counter $\theta$ are defined. When $\eta < \eta_{\min}$, $\theta$ is set to $\theta + 1$. If $\theta < \theta_{\max}$, the algorithm terminates. | Improvement factor $\eta$ and counter $\theta$ are defined. When $\eta < \eta_{\min}$, $\theta$ is set to $\theta + 1$. If $\theta < \theta_{\max}$, the algorithm terminates. |

## 4 Numerical experiments

The ability of the proposed MVSMCTS formulation is assessed by three weight minimization problems of truss structures. These problems include a 10-bar planar truss, a 25-bar spatial truss, and a 220-bar transmission tower. 10-bar planar truss is solved considering continuous sizing variable. The other two examples are solved for both sizing and shape variable. The truss structures are analyzed using the direct stiffness method. The computations entailed by the MVSMCTS formulation and the direct stiffness method are performed in the Python software using Intel Core i7 2.30 GHz processor and 40 GB memory.

### 4.1 Problem statement

The first example is the 10-bar planar truss shown in Fig. 11. This structure has been previously optimized by many researchers using a variety of techniques: harmony search (HS) (Lee and Geem 2004), PSO (Li et al. 2007), artificial bee colony (ABC) (Sonmez 2011), water evaporation optimization (WEO) (Kaveh and Bakhshpoori 2016), and political optimizer (PO) (Awad 2021).

The second test problem considers sizing and shape optimization of the 25-bar spatial truss shown in Fig. 12. This truss structure is previously designed by other methods such as GA (Wu and Chow 1995), PSO (Gholizadeh 2013), differential evolution (DE) (Ho-Huu et al. 2015), ABC (Jawad et al. 2021), and medalist learning algorithm (MLA) (He and Cui 2023).

The last example is the 220-bar transmission tower shown in Fig. 13. Various design parameters considered here are shown in Tables 2-9.

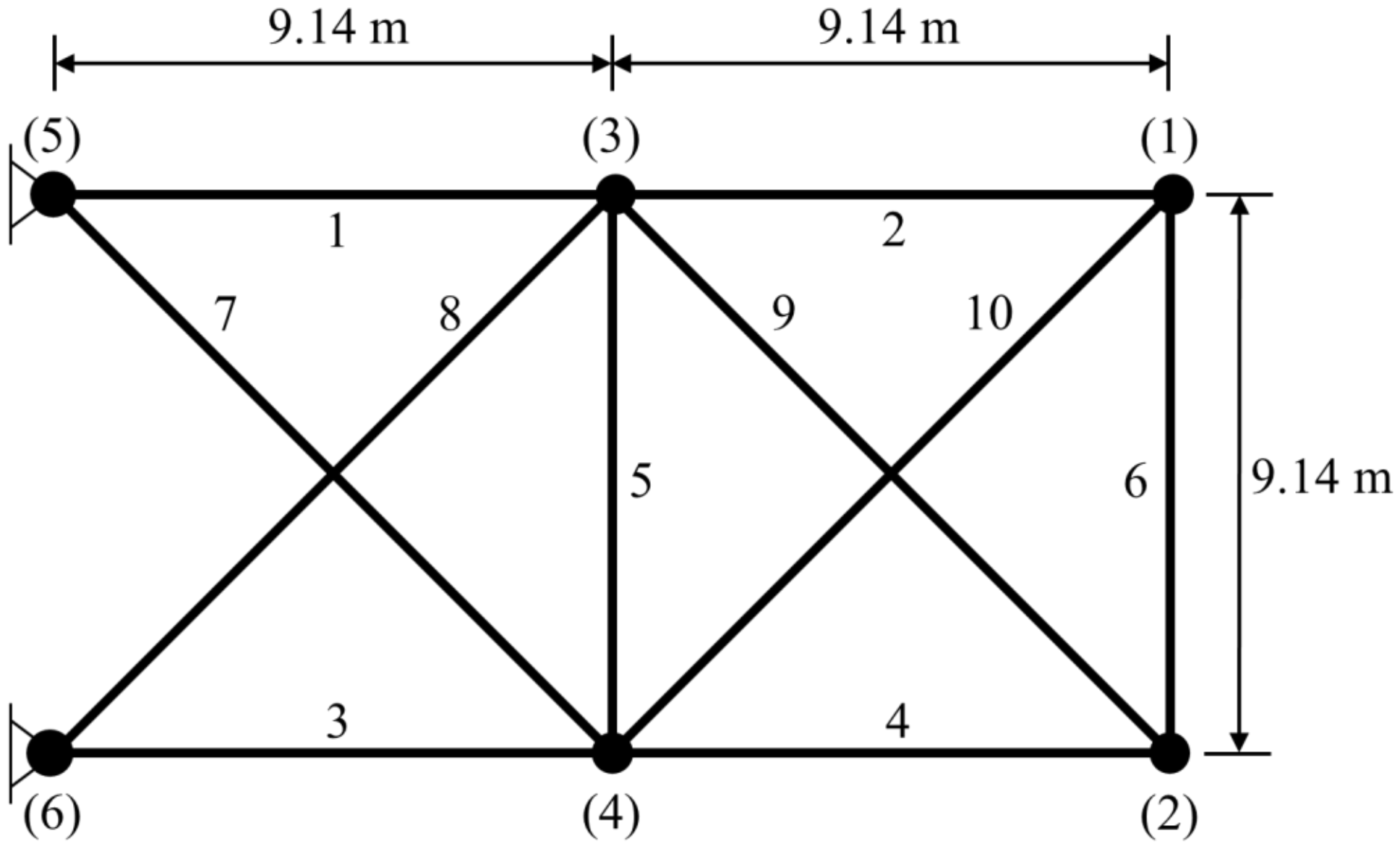

**Fig. 11** Schematic of the 10-bar planar truss structure

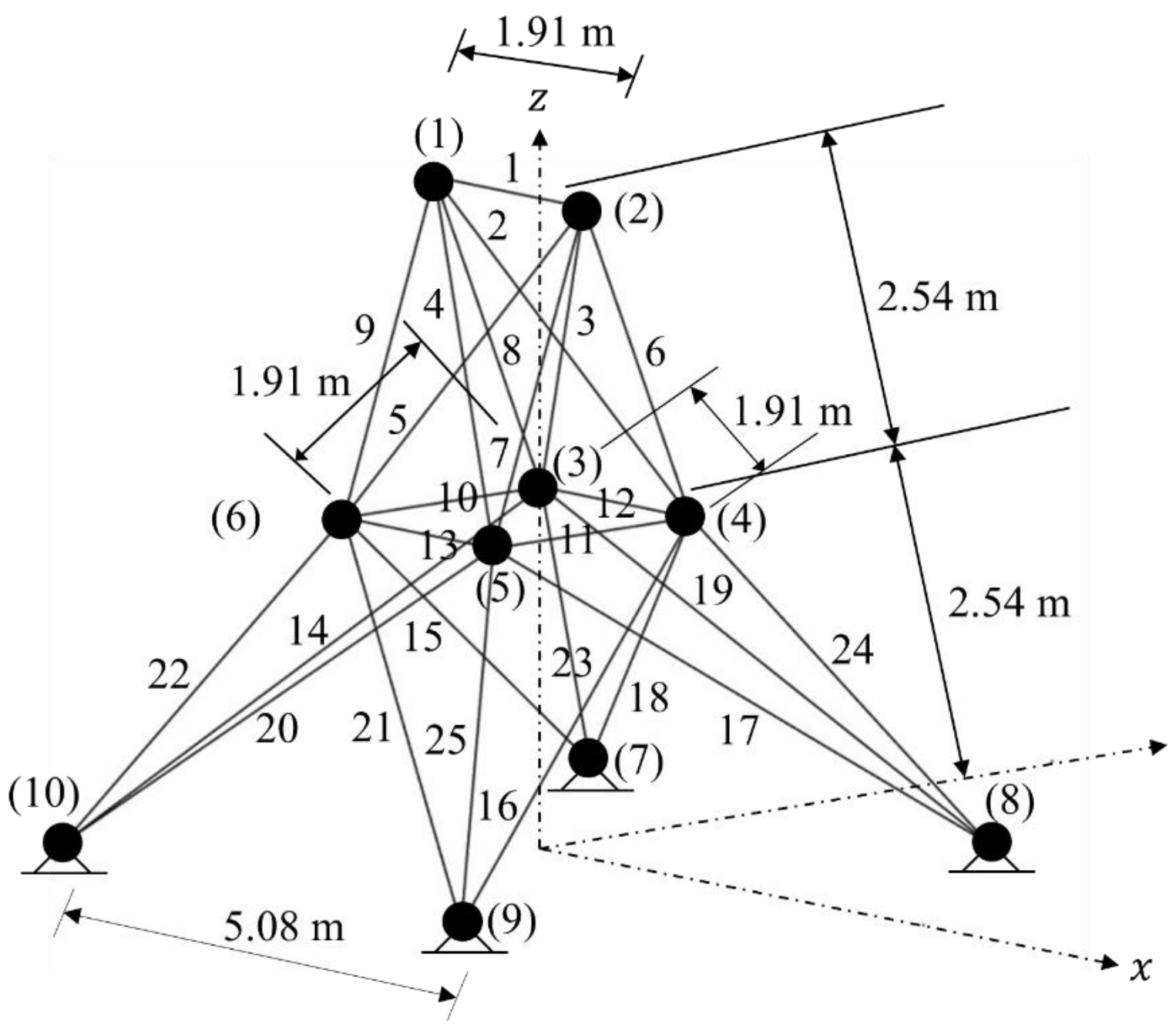

**Fig. 12** Schematic of the 25-bar spatial truss structure

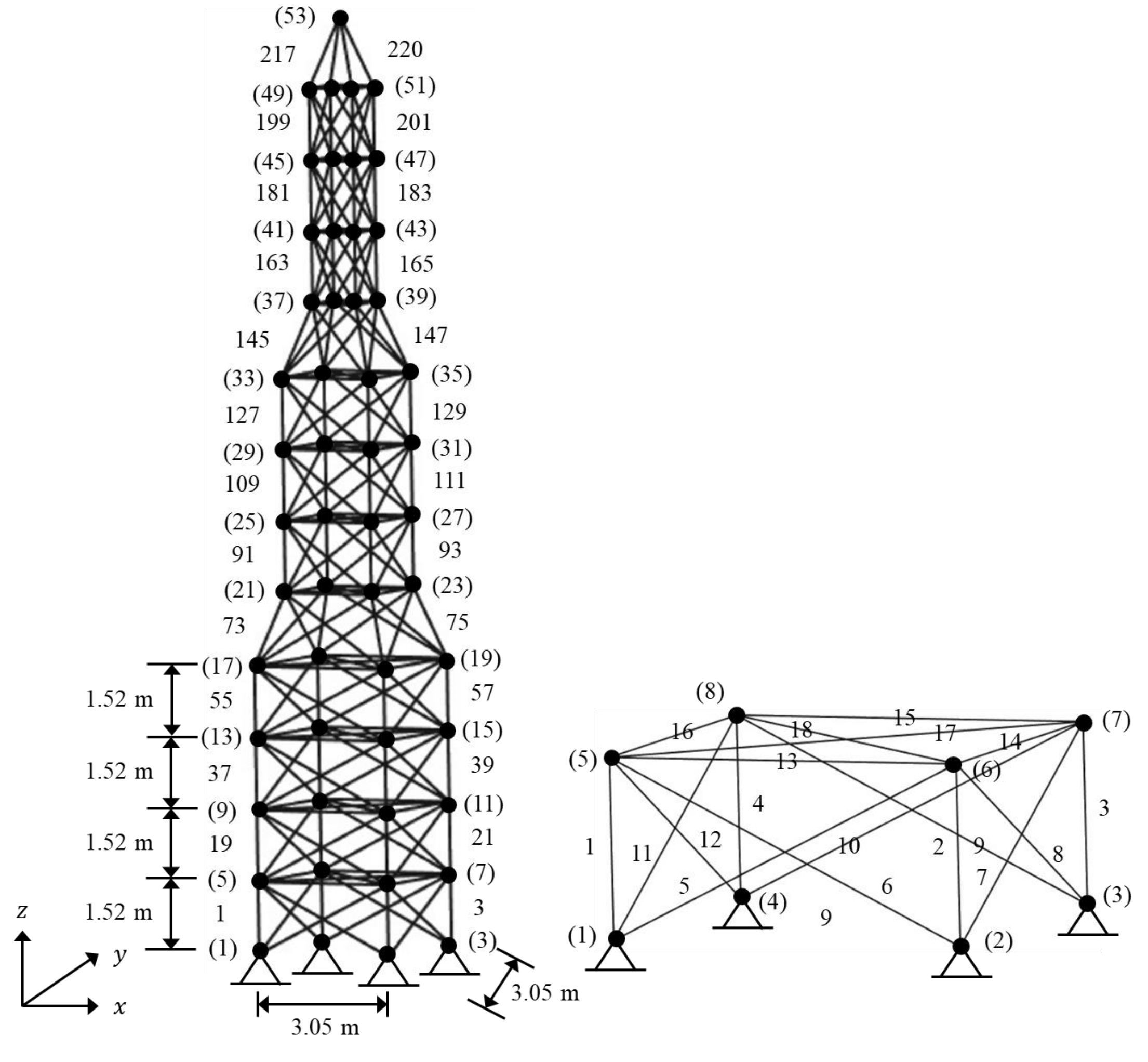

**Fig. 13** Schematic of the 220-bar transmission tower

**Table 2** Design consideration of the truss optimization problems

| | 10-bar planar truss | 25-bar spatial truss | 220-bar transmission tower |
|---|---|---|---|
| Design variables | $X_{i_X}, i_X = 1, 2, \dots, 10$ | $X_{i_X}, i_X = 1, 2, \dots, 8$<br>$Y_{i_Y}, i_Y = 1, 2, \dots, 5$ | $X_{i_X}, i_X = 1, 2, \dots, 49$<br>$Y_{i_Y}, i_Y = 1, 2, \dots, 21$ |
| Material density (kg/m$^3$) | 2767.99 | 2767.99 | 7860.00 |
| Modulus of elasticity (GPa) | 68.95 | 68.95 | 207.00 |
| Stress limitation (MPa) | ±172.37 | ±275.79 | ±180.00 |
| Displacement limitation (mm) | ±50.80 | ±8.89 | ±6.35 |
| Load (kN) | $P_2^y = P_4^y = -444.82$ | Table 4 | Table 5 |

**Table 3** Member grouping of the truss optimization problems

| Truss problems | Member group |
|---|---|
| 10-bar planar truss | No member group |
| 25-bar spatial truss | 25 truss members are grouped into 8 design variables, as follows: (1) 1, (2) 2–5, (3) 6–9, (4) 10–11, (5) 12–13, (6) 14–17, (7) 18–21, and (8) 22–25 |
| 220-bar transmission tower | Table 8 |

**Table 4** Load case for the 25-bar spatial truss structure

| Nodes $k$ | Load case (kN) | | |
|---|---|---|---|
| | $P_k^x$ | $P_k^y$ | $P_k^z$ |
| 1 | 4.45 | −44.48 | −44.48 |
| 2 | 0.00 | −44.48 | −44.48 |
| 3 | 2.22 | 0.00 | 0.00 |
| 6 | 2.67 | 0.00 | 0.00 |

**Table 5** Load case for the 220-bar transmission tower

| Nodes $k$ | Load case 1 (kN) | | | Load case 2 (kN) | | |
|---|---|---|---|---|---|---|
| | $P_k^x$ | $P_k^y$ | $P_k^z$ | $P_k^x$ | $P_k^y$ | $P_k^z$ |
| 49 | 0.00 | 0.00 | 0.00 | 0.00 | 0.00 | −400.00 |
| 50 | 0.00 | 0.00 | 0.00 | 0.00 | 0.00 | −400.00 |
| 51 | 0.00 | 0.00 | 0.00 | 0.00 | 0.00 | −400.00 |
| 52 | 0.00 | 0.00 | 0.00 | 0.00 | 0.00 | −400.00 |
| 53 | 0.00 | 0.00 | −2000.00 | 0.00 | 0.00 | −400.00 |

**Table 6** The lower and upper limit for sizing and shape

| | 10-bar planar truss | 25-bar spatial truss |
|---|---|---|
| $X_{i_X,\min}$ (mm$^2$) | 64.52 | 64.52 |
| $X_{i_X,\max}$ (mm$^2$) | 22582.00 | 2193.68 |
| $Y_{i_Y,\min}$ (m)<br>$Y_{i_Y,\max}$ (m) | N/A | $0.51 \le x_4 = x_5 = -x_3 = -x_6 \le 1.52$,<br>$1.02 \le y_3 = y_4 = -y_5 = -y_6 \le 2.03$,<br>$2.29 \le z_3 = z_4 = z_5 = z_6 \le 3.30$,<br>$1.02 \le x_8 = x_9 = -x_7 = -x_{10} \le 2.03$,<br>$2.54 \le y_7 = y_8 = -y_9 = -y_{10} \le 3.56$ |

| | 220-bar transmission tower |
|---|---|
| $X_{i_X,\min}$ (mm$^2$) | 71.61 |
| $X_{i_X,\max}$ (mm$^2$) | 21612.86 |
| $Y_{i_Y,\min}$ (m)<br>$Y_{i_Y,\max}$ (m) | $-0.25 \le x_5 = y_5 = y_6 = x_8 = x_9 = y_9 = y_{10} = x_{12} = x_{13} = y_{13} = y_{14} = x_{16} = x_{17} = y_{17} = y_{18} = x_{20} \le 0.25$,<br>$2.79 \le x_6 = x_7 = y_7 = y_8 = x_{10} = x_{11} = y_{11} = y_{12} = x_{14} = x_{15} = y_{15} = y_{16} = x_{18} = x_{19} = y_{19} = y_{20} \le 3.30$,<br>$0.25 \le x_{21} = y_{21} = y_{22} = x_{24} = x_{25} = y_{25} = y_{26} = x_{28} = x_{29} = y_{29} = y_{30} = x_{32} = x_{33} = y_{33} = y_{34} = x_{36} \le 0.76$,<br>$2.29 \le x_{22} = x_{23} = y_{23} = y_{24} = x_{26} = x_{27} = y_{27} = y_{28} = x_{30} = x_{31} = y_{31} = y_{32} = x_{34} = x_{35} = y_{35} = y_{36} \le 2.79$,<br>$0.76 \le x_{37} = y_{37} = y_{38} = x_{40} = x_{41} = y_{41} = y_{42} = x_{44} = x_{45} = y_{45} = y_{46} = x_{48} = x_{49} = y_{49} = y_{50} = x_{52} \le 1.27$,<br>$1.78 \le x_{38} = x_{39} = y_{39} = y_{40} = x_{42} = x_{43} = y_{43} = y_{44} = x_{46} = x_{47} = y_{47} = y_{48} = x_{50} = x_{51} = y_{51} = y_{52} \le 2.29$,<br>$1.27 \le x_{53} \le 1.78$, $1.27 \le y_{53} \le 1.78$, $1.27 \le z_5 = z_6 = z_7 = z_8 \le 1.78$, $2.79 \le z_9 = z_{10} = z_{11} = z_{12} \le 3.30$,<br>$4.32 \le z_{13} = z_{14} = z_{15} = z_{16} \le 4.83$, $5.84 \le z_{17} = z_{18} = z_{19} = z_{20} \le 6.35$, $7.34 \le z_{21} = z_{22} = z_{23} = z_{24} \le 7.87$,<br>$8.89 \le z_{25} = z_{26} = z_{27} = z_{28} \le 9.40$, $10.41 \le z_{29} = z_{30} = z_{31} = z_{32} \le 10.92$, $11.94 \le z_{33} = z_{34} = z_{35} = z_{36} \le 12.45$,<br>$13.46 \le z_{37} = z_{38} = z_{39} = z_{40} \le 13.97$, $14.99 \le z_{41} = z_{42} = z_{43} = z_{44} \le 15.49$, $16.51 \le z_{45} = z_{46} = z_{47} = z_{48} \le 17.02$,<br>$18.03 \le z_{49} = z_{50} = z_{51} = z_{52} \le 18.54$, $19.56 \le z_{53} \le 20.07$ |

**Table 7** Discrete set for sizing and shape

| Truss problems | Discrete set |
|---|---|
| 25-bar spatial truss | $\mathbf{D}_{X,i_X}$ ={0.65, 1.30, 1.95,…, 15.60, 16.25, 16.90,…,22.10} cm$^2$, $i_X = 1, 2, \ldots, 8$<br>$\mathbf{D}_{Y,1}$ ={0.51, 0.57, 0.64, 0.70,…,1.33, 1.40, 1.46, 1.52} m<br>$\mathbf{D}_{Y,2}$ ={1.02, 1.08, 1.14, 1.21,…,1.84, 1.91, 1.97, 2.03} m<br>$\mathbf{D}_{Y,3}$ ={2.29, 2.35, 2.41, 2.48,…,3.11, 3.18, 3.24, 3.30} m<br>$\mathbf{D}_{Y,4}$ ={1.02, 1.08, 1.14, 1.21,…,1.84, 1.91, 1.97, 2.03} m<br>$\mathbf{D}_{Y,5}$ ={2.54, 2.60, 2.67, 2.73,…,3.37, 3.43, 3.49, 3.56} m |
| 220-bar transmission tower | $\mathbf{D}_{X,i_X}$ ($i_X = 1, 2, \ldots, 49$) refers to Table 9.<br>$\mathbf{D}_{Y,1}$ ={−0.25, −0.22, −0.19, −0.16,…,0.16, 0.19, 0.22, 0.25} m<br>$\mathbf{D}_{Y,2}$ = $\mathbf{D}_{Y,10}$ ={2.79, 2.83, 2.86, 2.89,…,3.21, 3.24, 3.27, 3.30} m<br>$\mathbf{D}_{Y,3}$ ={0.25, 0.29, 0.32, 0.35,…,0.67, 0.70, 0.73, 0.76} m<br>$\mathbf{D}_{Y,4}$ ={2.29, 2.32, 2.35, 2.38,…,2.70, 2.73, 2.76, 2.79} m<br>$\mathbf{D}_{Y,5}$ ={0.76, 0.79, 0.83, 0.86,…,1.17, 1.21, 1.24, 1.27} m<br>$\mathbf{D}_{Y,6}$ ={1.78, 1.81, 1.84, 1.87,…,2.19, 2.22, 2.25, 2.29} m<br>$\mathbf{D}_{Y,7}$ = $\mathbf{D}_{Y,8}$ = $\mathbf{D}_{Y,9}$ ={1.27, 1.30, 1.33, 1.37,…,1.68, 1.71, 1.75, 1.78} m<br>$\mathbf{D}_{Y,11}$ ={4.32, 4.35, 4.38, 4.41,…,4.73, 4.76, 4.79, 4.83} m<br>$\mathbf{D}_{Y,12}$ ={5.84, 5.87, 5.91, 5.94,…,6.25, 6.29, 6.32, 6.35} m<br>$\mathbf{D}_{Y,13}$ ={7.37, 7.40, 7.43, 7.46,…,7.78, 7.81, 7.84, 7.87} m<br>$\mathbf{D}_{Y,14}$ ={8.89, 8.92, 8.95, 8.99,…,9.30, 9.33, 9.37, 9.40} m<br>$\mathbf{D}_{Y,15}$ ={10.41, 10.45, 10.48, 10.51,…,10.83, 10.86, 10.89, 10.92} m<br>$\mathbf{D}_{Y,16}$ ={11.94, 11.97, 12.00, 12.03,…,12.35, 12.38, 12.41, 12.45} m<br>$\mathbf{D}_{Y,17}$ ={13.53, 13.56, 13.59, 13.62,…,13.87, 13.91, 13.94, 13.97} m<br>$\mathbf{D}_{Y,18}$ ={14.99, 15.02, 15.05, 15.08,…,15.40, 15.43, 15.46, 15.49} m<br>$\mathbf{D}_{Y,19}$ ={16.51, 16.54, 16.57, 16.61,…,16.92, 16.95, 16.99, 17.02} m<br>$\mathbf{D}_{Y,20}$ ={18.07, 18.10, 18.13, 18.16,…,18.45, 18.48, 18.51, 18.54} m<br>$\mathbf{D}_{Y,21}$ ={19.59, 19.62, 19.65, 19.69,…,19.97, 20.00, 20.03, 20.07} m |

**Table 8** Member grouping of the 220-bar transmission tower

| Number | Members | Number | Members | Number | Members | Number | Members |
|---|---|---|---|---|---|---|---|
| 1 | 1–4 | 14 | 59–66 | 27 | 121–124 | 40 | 179–180 |
| 2 | 5–12 | 15 | 67–70 | 28 | 125–126 | 41 | 181–184 |
| 3 | 13–16 | 16 | 71–72 | 29 | 127–130 | 42 | 185–192 |
| 4 | 17–18 | 17 | 73–76 | 30 | 131–138 | 43 | 193–196 |
| 5 | 19–22 | 18 | 77–84 | 31 | 139–142 | 44 | 197–198 |
| 6 | 23–30 | 19 | 85–88 | 32 | 143–144 | 45 | 199–202 |
| 7 | 31–34 | 20 | 89–90 | 33 | 145–148 | 46 | 203–210 |
| 8 | 35–36 | 21 | 91–94 | 34 | 149–156 | 47 | 211–214 |
| 9 | 37–40 | 22 | 95–102 | 35 | 157–160 | 48 | 215–216 |
| 10 | 41–48 | 23 | 103–106 | 36 | 161–162 | 49 | 217–220 |
| 11 | 49–52 | 24 | 107–108 | 37 | 163–166 | | |
| 12 | 53–54 | 25 | 109–112 | 38 | 167–174 | | |
| 13 | 55–58 | 26 | 113–120 | 39 | 175–178 | | |

**Table 9** Discrete values available for cross-sectional areas from AISC norm

| Number | Area (mm$^2$) | Number | Area (mm$^2$) | Number | Area (mm$^2$) | Number | Area (mm$^2$) |
|---|---|---|---|---|---|---|---|
| 1 | 71.61 | 17 | 1008.39 | 33 | 2477.41 | 49 | 7419.34 |
| 2 | 90.97 | 18 | 1045.16 | 34 | 2496.77 | 50 | 8709.66 |
| 3 | 126.45 | 19 | 1161.29 | 35 | 2503.22 | 51 | 8967.72 |
| 4 | 161.29 | 20 | 1283.87 | 36 | 2696.77 | 52 | 9161.27 |
| 5 | 198.06 | 21 | 1374.19 | 37 | 2722.58 | 53 | 9999.98 |
| 6 | 252.26 | 22 | 1535.48 | 38 | 2896.77 | 54 | 10322.56 |
| 7 | 285.16 | 23 | 1690.32 | 39 | 2961.28 | 55 | 10903.20 |
| 8 | 363.23 | 24 | 1696.77 | 40 | 3096.77 | 56 | 12129.01 |
| 9 | 388.39 | 25 | 1858.06 | 41 | 3206.45 | 57 | 12838.68 |
| 10 | 494.19 | 26 | 1890.32 | 42 | 3303.22 | 58 | 14193.52 |
| 11 | 506.45 | 27 | 1993.54 | 43 | 3703.22 | 59 | 14774.16 |
| 12 | 641.29 | 28 | 2019.35 | 44 | 4658.06 | 60 | 15806.42 |
| 13 | 645.16 | 29 | 2180.64 | 45 | 5141.93 | 61 | 17096.74 |
| 14 | 792.26 | 30 | 2238.71 | 46 | 5503.22 | 62 | 18064.48 |
| 15 | 816.77 | 31 | 2290.32 | 47 | 5999.99 | 63 | 19354.80 |
| 16 | 940.00 | 32 | 2341.93 | 48 | 6999.99 | 64 | 21612.86 |

### 4.2 The effect of expansion of search region on the optimal solution for continuous sizing and shape variable

Figs. 14 and 15 illustrate the convergence histories for the 10-bar planar truss and the 25-bar spatial truss under different constants $A_X$ and $A_Y$. It is worth mentioning that all design variables are continuous. The figures demonstrate that the proposed algorithm obtains the best solution when $A_X$ and $A_Y$ are all equal to 0.5. It indicates that the search region in each round must be expanded from the center.

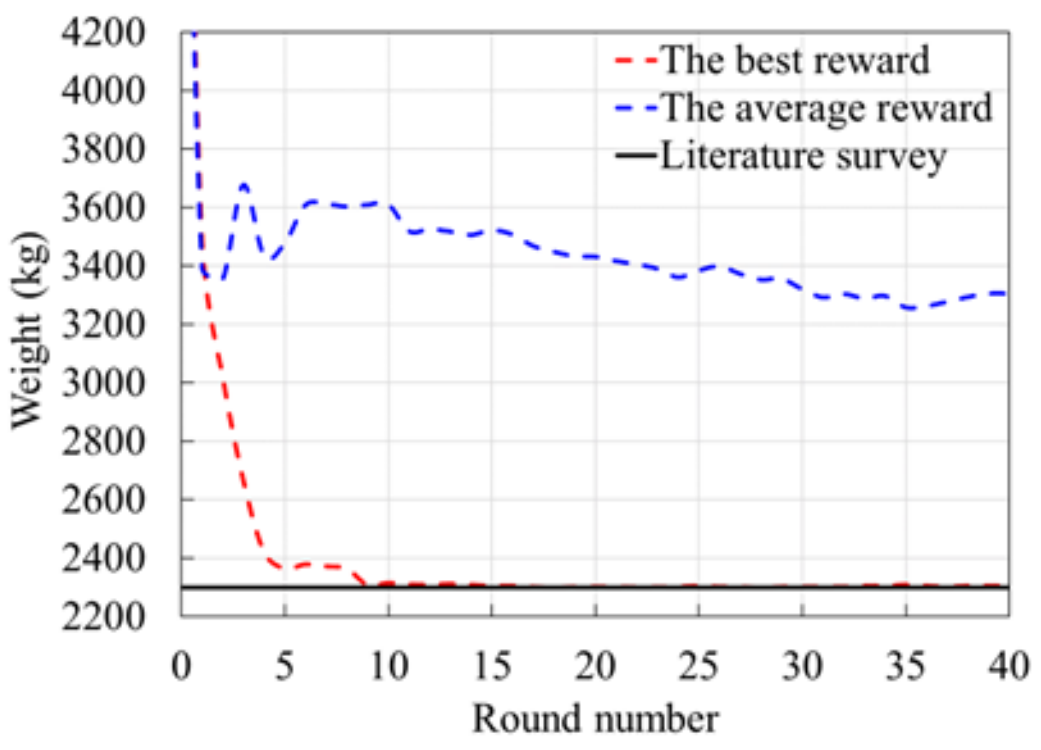

**Fig. 14** Convergence histories for the 10-bar planar truss considering continuous sizing variable under different constant parameter $A_X$

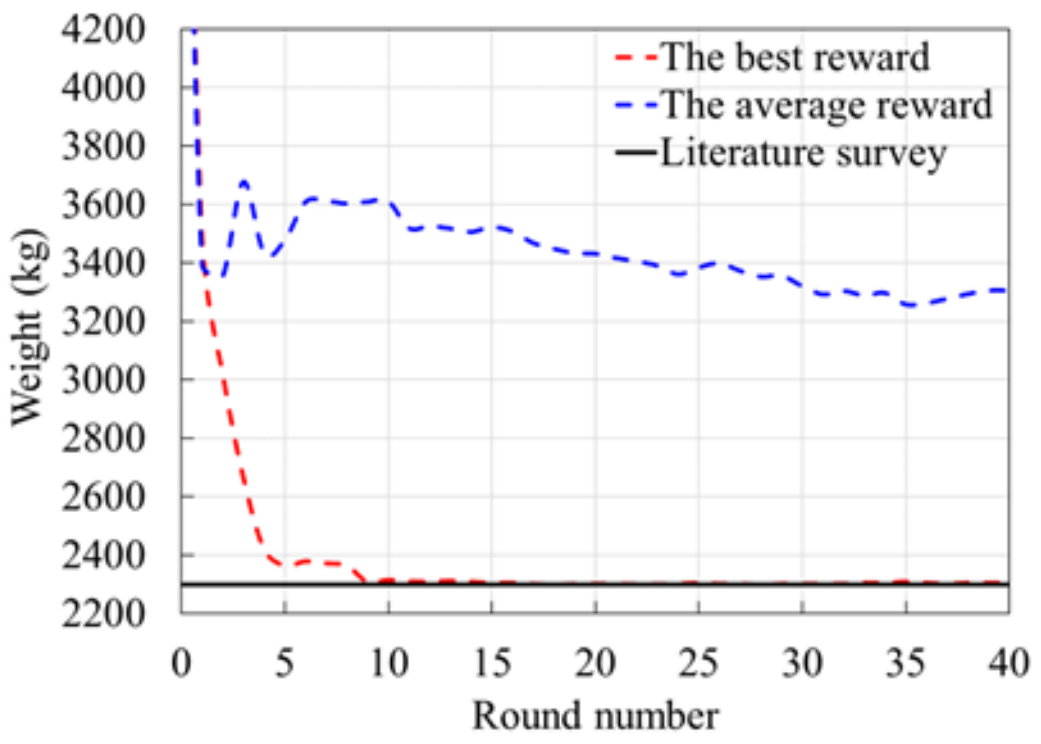

**Fig. 15** Convergence histories for the 25-bar spatial truss considering continuous sizing and shape variable under different constant parameters $A_X$ and $A_Y$

### 4.3 The effect of different types of accelerating techniques for continuous variable on the computation time

Table 10 compares the efficiency of IMCTS without and with different types of accelerating techniques based on CPU time. It can be seen that geometric decay is 1.5, 2.2, and 4.8 times faster than linear decrease, step reduction, and without accelerating technique.

**Table 10.** CPU time of IMCTS formulation without and with different types of accelerating techniques

| | Without accelerating technique | Geometric decay | Linear decrease | Step reduction |
|---|---|---|---|---|
| 10-bar planar truss (Case 1) | 223.49 sec | 46.53 sec | 62.94 sec | 102.47 sec |
| 10-bar planar truss (Case 2) | 380.13 sec | 86.02 sec | 102.45 sec | 264.03 sec |
| 72-bar spatial truss (Case 1) | 12749.90 sec | 2549.98 sec | 3470.55 sec | 4978.25 sec |
| 72-bar spatial truss (Case 2) | 36154.62 sec | 8408.05 sec | 14851.48 sec | 20872.75 sec |

### 4.4. Investigation of convergence history: 10-bar planar truss

Figure 18 illustrates the comparison of convergence histories for the 10-bar planar truss under the best and average reward. Figure 18(a) and (b) for Case 1 and Case 2 demonstrates that the proposed method obtains the best solution at 13 and 15 rounds under the best reward. However, this method does not detect the best solution after 40 rounds under the average reward.

The comparison of convergence histories for the 10-bar planar truss under different parameter $\alpha$ is shown in Figure 19. From Figure 19(a) and (b), this algorithm obtains the best solution at 13 and 18 rounds under minimum weight for Case 1 and Case 2. However, for Case 1 and Case 2, this algorithm does not detect the best solution after 40 rounds under maximum weight.

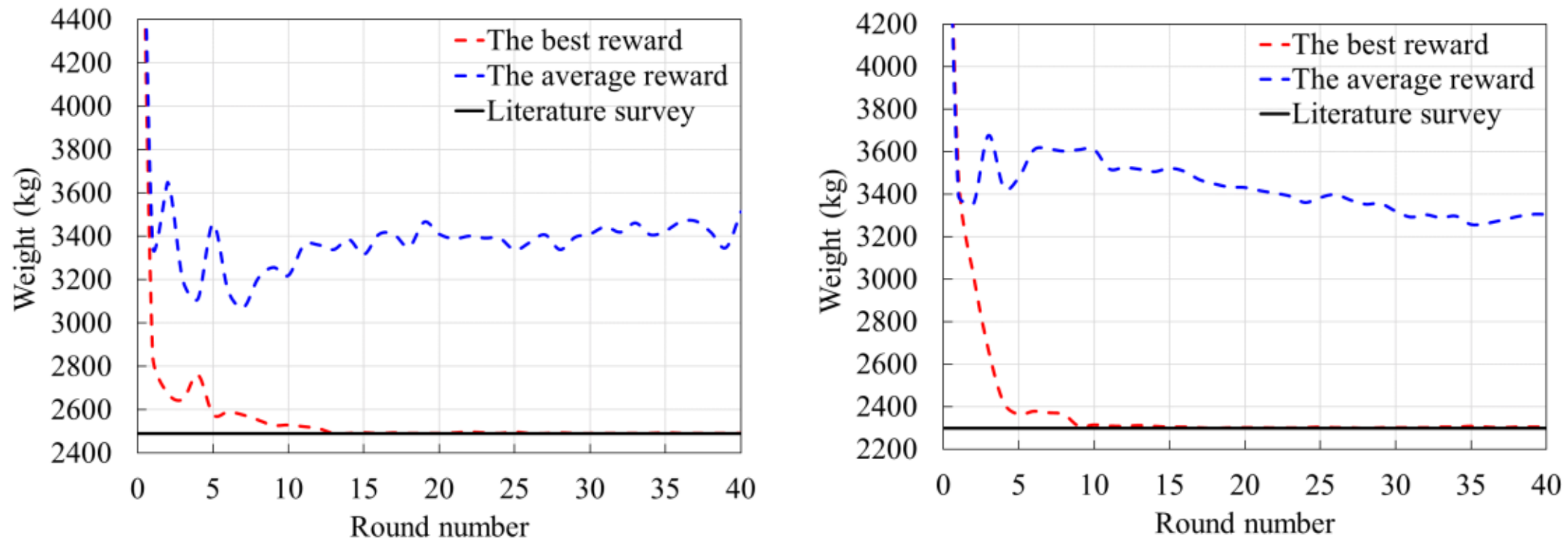

(a) (b)

**Figure 18.** Comparison of the convergence histories for 10-bar planar truss under the best and average reward for (a) Case 1 and (b) Case 2.

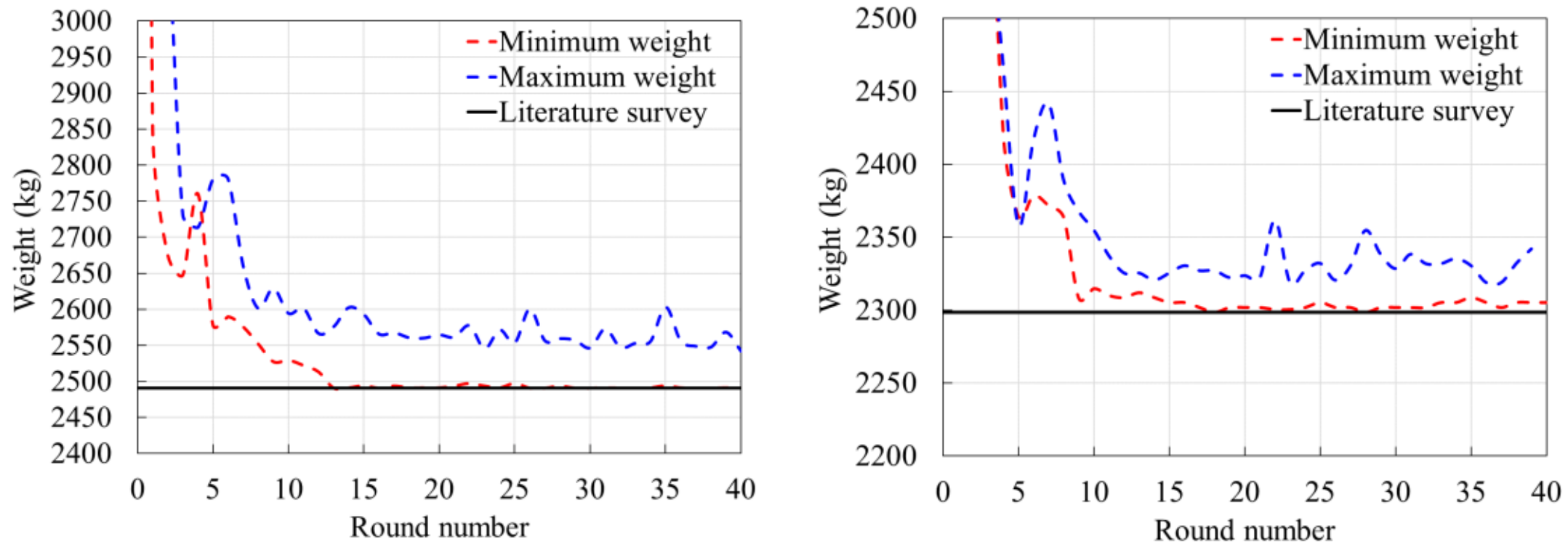

(a) (b)

**Figure 19.** Comparison of the convergence histories for 10-bar planar truss under parameter $\alpha$ equal to minimum and maximum weight for (a) Case 1 and (b) Case 2.

## 4.5 Solution accuracy

Three truss optimization problems, including a 15-bar planar truss, a 25-bar spatial truss, and a 220-bar transmission tower are used to validate the solution accuracy by comparing the results which have been previously investigated by other researchers.

Tables 9-10 demonstrate the comparison of optimal design results for a 15-bar planar truss and a 25-bar spatial truss. It is shown that the IMCTS formulation outperforms all metaheuristic algorithms in terms of the lightest weight. Optimal layouts after the optimization with the IMCTS formulation are given in Figs. 15-17. Figs. 18-20 indicate that the optimal solution assessed by the IMCTS formulation does not violate normal stresses and nodal displacement constraints.

**Table 10** Comparison of optimal designs for the 25-bar spatial truss structure

| Design variables | GA | PSO | FA | DE | ABC | UMCTS |
|---|---|---|---|---|---|---|
| Sizing optimization/Cross-sectional area ($mm^2$) | | | | | | |
| 1 | 0.65 | 0.65 | 0.65 | 0.65 | 0.65 | 0.65 |
| 2-5 | 1.30 | 0.65 | 0.65 | 0.65 | 0.65 | 0.65 |
| 6-9 | 7.15 | 7.15 | 5.85 | 5.85 | 6.45 | 6.50 |
| 10-11 | 1.30 | 0.65 | 0.65 | 0.65 | 0.65 | 0.65 |
| 12-13 | 1.95 | 2.60 | 0.65 | 0.65 | 0.65 | 0.65 |
| 14-17 | 0.65 | 0.65 | 0.65 | 0.65 | 0.65 | 0.65 |
| 18-21 | 1.30 | 2.60 | 0.65 | 0.65 | 0.65 | 0.65 |
| 22-25 | 5.85 | 4.55 | 6.45 | 6.45 | 5.85 | 5.85 |
| Shape optimization/Nodal coordinate (mm) | | | | | | |
| $x_4$ | 1.04 | 0.70 | 0.95 | 0.94 | 0.92 | 0.95 |
| $y_3$ | 1.36 | 1.31 | 1.42 | 1.49 | 1.39 | 1.39 |
| $z_3$ | 3.16 | 3.30 | 3.22 | 3.12 | 3.30 | 3.29 |
| $x_8$ | 1.29 | 1.08 | 1.27 | 1.25 | 1.32 | 1.30 |
| $y_7$ | 3.34 | 3.37 | 3.46 | 3.47 | 3.56 | 3.49 |
| Weight (kg) | 61.78 | 58.61 | 53.90 | 53.87 | 53.22 | 52.96 |

**Table 11** Optimal design for 220-bar transmission tower

| Sizing optimization | | | | | | | |
|---|---|---|---|---|---|---|---|
| Number | Cross-sectional area ($mm^2$) | Number | Cross-sectional area ($mm^2$) | Number | Cross-sectional area ($mm^2$) | Number | Cross-sectional area ($mm^2$) |
| 1 | 8709.66 | 14 | 71.61 | 27 | 71.61 | 40 | 90.97 |
| 2 | 71.61 | 15 | 198.06 | 28 | 71.61 | 41 | 7419.34 |
| 3 | 90.97 | 16 | 1858.06 | 29 | 8709.66 | 42 | 71.61 |
| 4 | 1690.32 | 17 | 8967.72 | 30 | 71.61 | 43 | 71.61 |
| 5 | 8967.72 | 18 | 71.61 | 31 | 641.29 | 44 | 71.61 |
| 6 | 71.61 | 19 | 1161.29 | 32 | 2180.64 | 45 | 7419.34 |
| 7 | 90.97 | 20 | 645.16 | 33 | 8709.66 | 46 | 71.61 |
| 8 | 71.61 | 21 | 8709.66 | 34 | 90.97 | 47 | 1690.32 |
| 9 | 8709.66 | 22 | 71.61 | 35 | 494.19 | 48 | 252.26 |
| 10 | 90.97 | 23 | 71.61 | 36 | 2503.22 | 49 | 8709.66 |
| 11 | 90.97 | 24 | 71.61 | 37 | 8709.66 | | |
| 12 | 161.29 | 25 | 8967.72 | 38 | 71.61 | | |
| 13 | 9161.27 | 26 | 71.61 | 39 | 90.97 | | |
| Shape optimization | | | | | | | |
| Number | Nodal coordinate (mm) | Number | Nodal coordinate (mm) | Number | Nodal coordinate (mm) | Number | Nodal coordinate (mm) |
| 1 | 242.09 | 7 | 1492.16 | 13 | 7866.06 | 19 | 16545.712 |
| 2 | 2797.97 | 8 | 1491.10 | 14 | 9178.75 | 20 | 18085.15 |
| 3 | 612.50 | 9 | 1695.93 | 15 | 10414.00 | 21 | 19563.95 |
| 4 | 2399.00 | 10 | 2951.43 | 16 | 11938.00 | | |
| 5 | 1137.09 | 11 | 4353.72 | 17 | 13966.03 | | |
| 6 | 1861.70 | 12 | 5842.00 | 18 | 15387.90 | | |

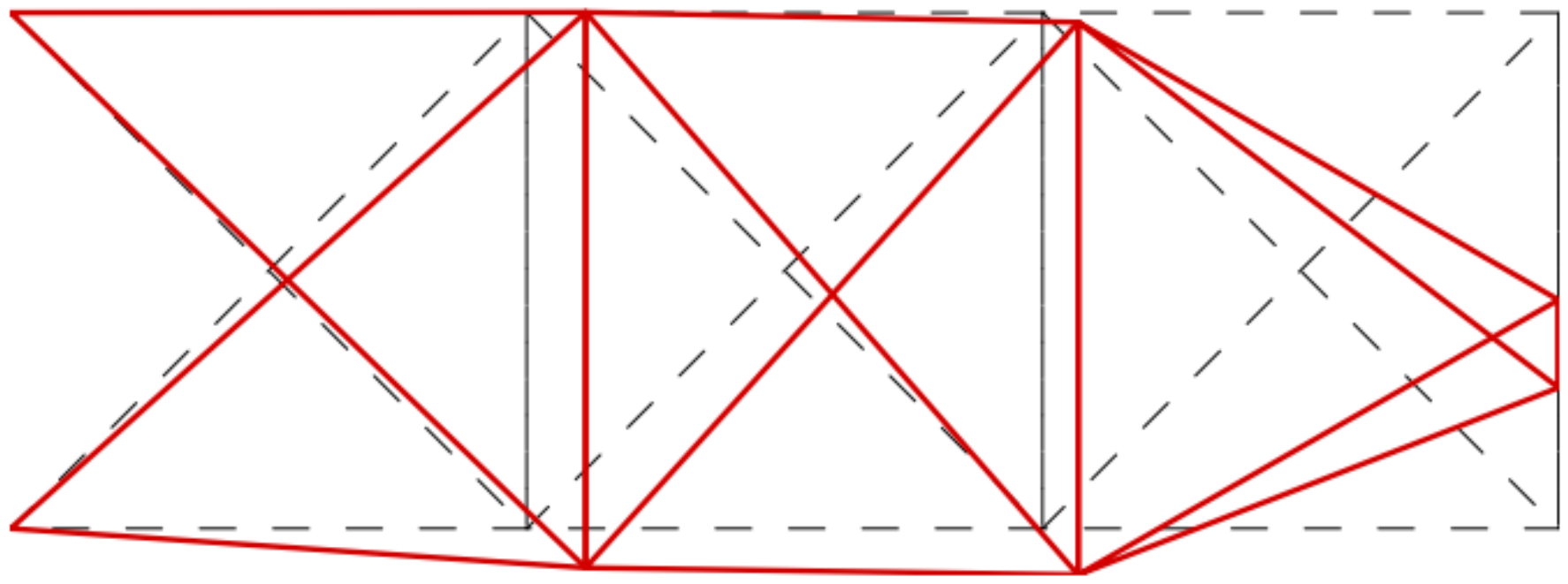

**Fig. 15** 15-bar planar truss structure: comparison of the optimized layout with the initial configuration of the truss

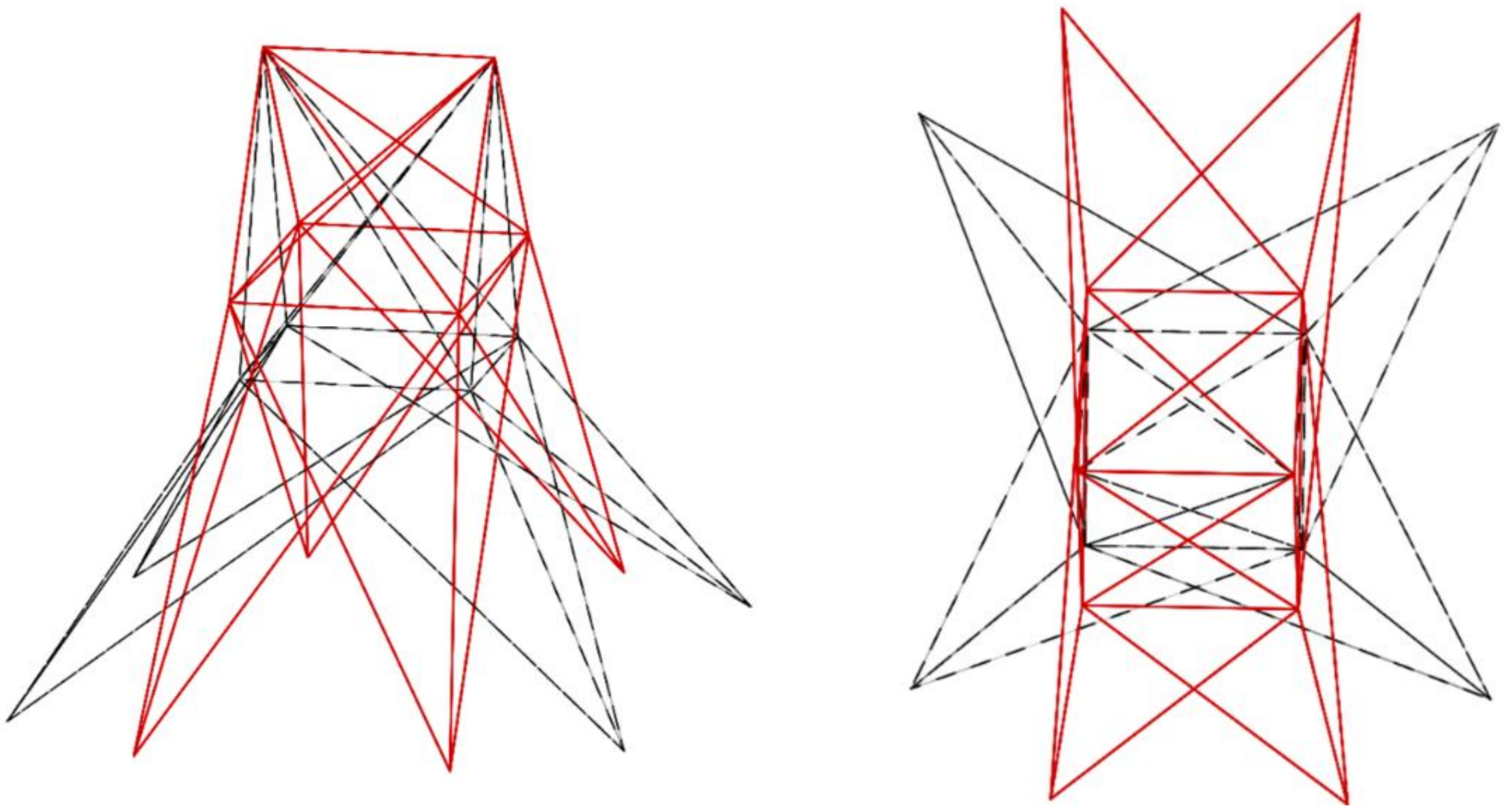

**Fig. 16** 25-bar spatial truss structure: comparison of the optimized layout with the initial configuration of the truss

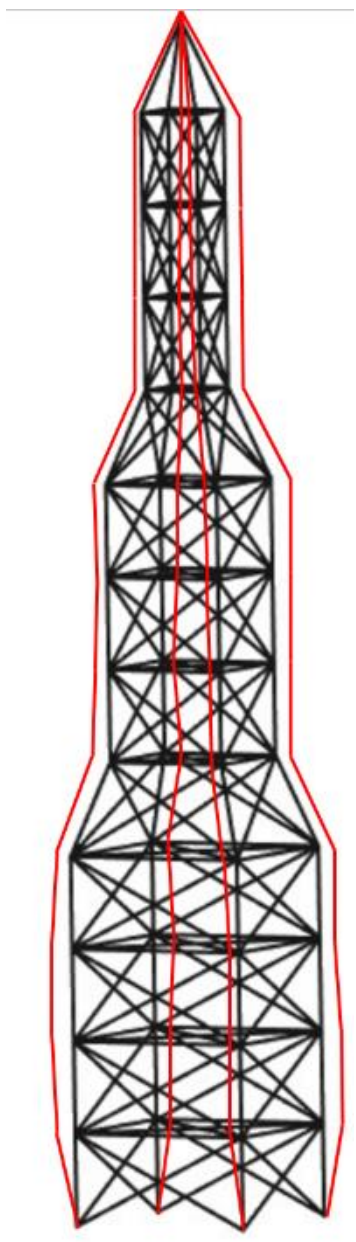

**Fig. 17** 220-bar transmission tower: comparison of the optimized layout with the initial configuration of the truss

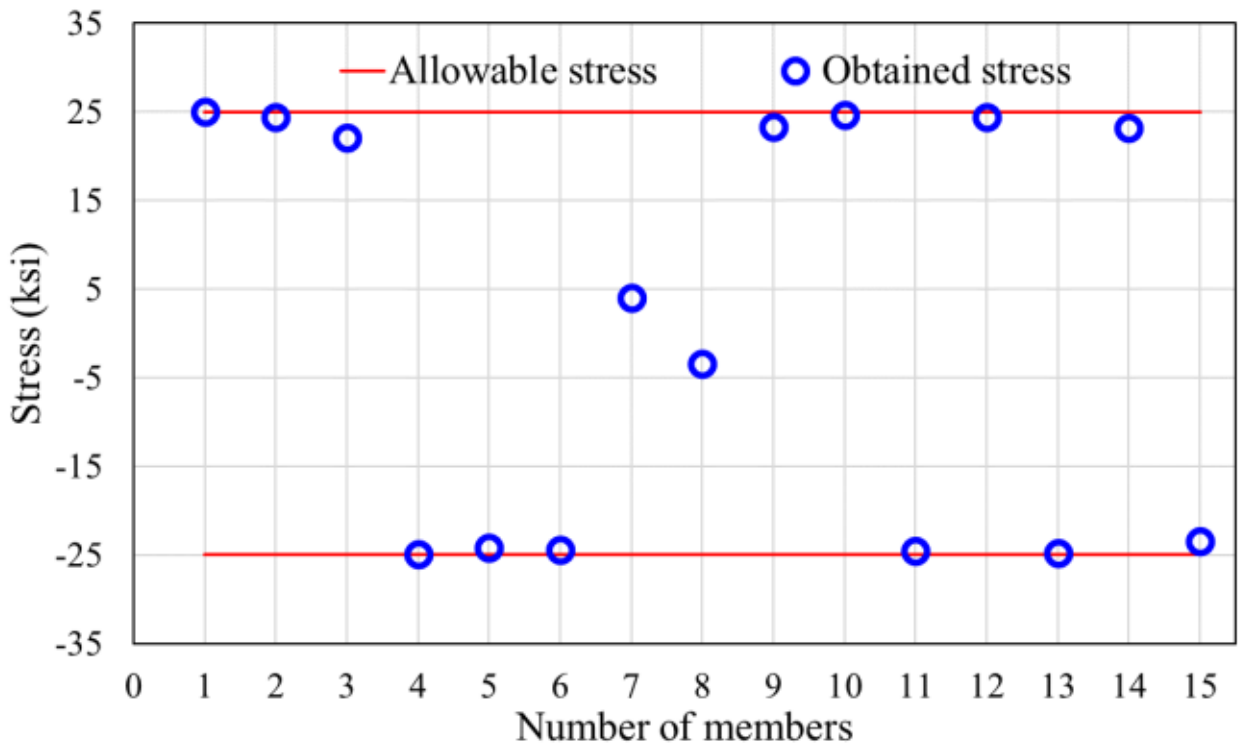

**Fig. 18** Constraints boundaries assessed at the optimal design of 15-bar planar truss structure by the IMCTS formulation

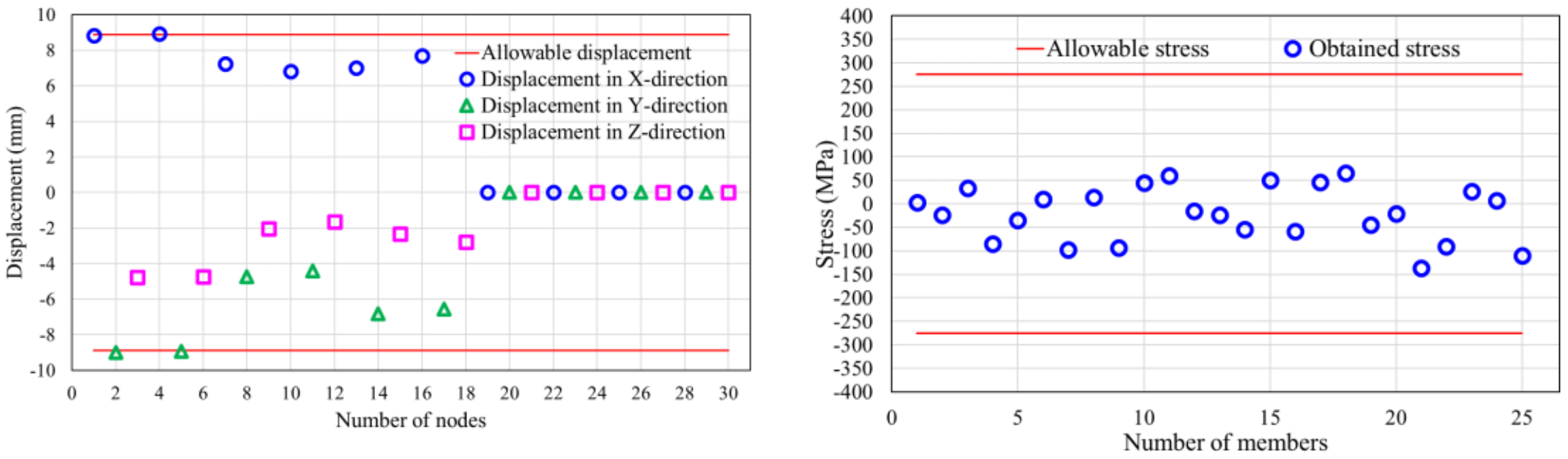

(a) (b)

**Fig. 19** Constraints boundaries assessed at the optimal design of 25-bar spatial truss structure by the IMCTS formulation for (a) displacement constraints (b) stress constraints

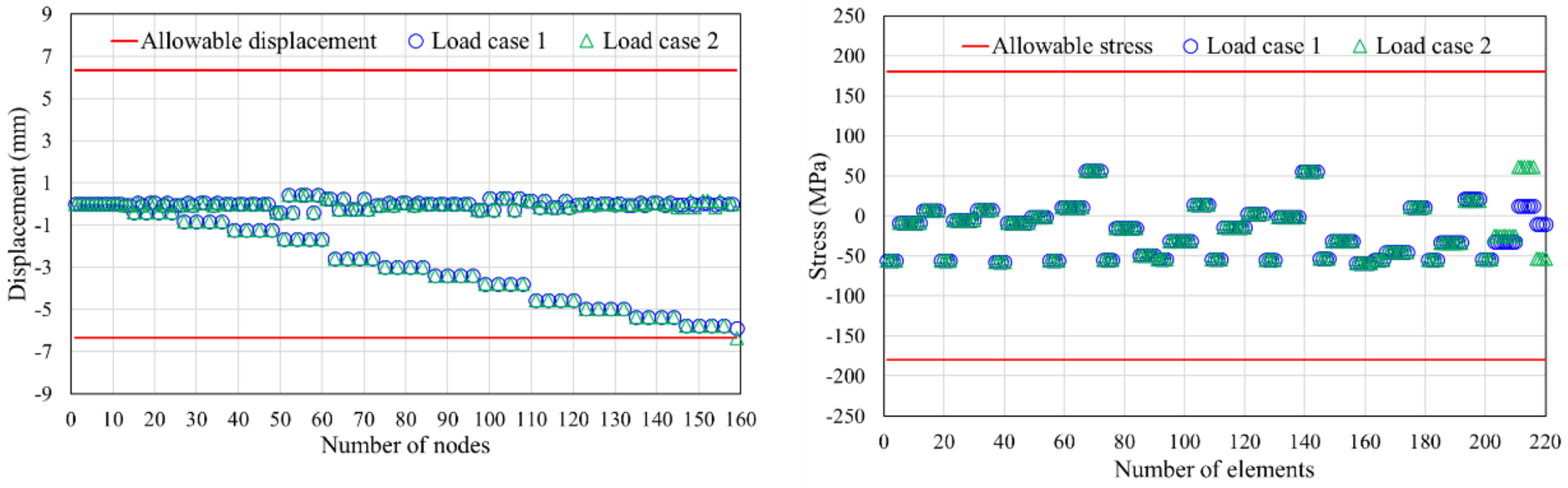

(a) (b)

**Fig. 20** Constraints boundaries assessed at the optimal design of 220-bar transmission tower by the IMCTS formulation for (a) displacement constraints (b) stress constraints

## 4.4 Solution stability

The statistical results for 15-bar planar truss, 25-bar spatial truss, and 220-bar transmission tower are obtained through 10 independent sampling to test the stability of this algorithm. The results are shown in Table 12, including the best, the worst, average, and standard deviation.

**Table 12** Statistical results of the investigated example

| Investigated example | Best weight | Worst weight | Average weight | Standard deviation |
|---|---|---|---|---|
| 15-bar planar truss | 34.40 kg | 36.39 kg | 35.66 kg | 0.59 |
| 25-bar spatial truss | 53.16 kg | 56.91 kg | 54.24 kg | 1.02 |
| 220-bar transmission tower | 6113.64 kg | 6302.15 kg | 6226.99 kg | 76.98 |

## 5 Summary and conclusions

In this paper, an RL algorithm using MVSMCTS formulation is developed to solve optimization problems considering various types of variables with single and mixed system. MVSMCTS formulation is based on IMCTS formulation for discrete sizing optimization. Sizing and shape optimization of truss structures is utilized as an example for mixed variable structural optimization. For this problem, the design variables are the cross-sectional areas of the members and the nodal coordinates of the joints. It is difficult to find an optimal solution in a reasonable time because various types of variables are of fundamentally different nature. The following conclusions are obtained:

(1) The computational framework of update process for discrete variable is also suitable for continuous variable. In order to adopt update process in continuous space, uniform meshes are generated automatically in search region.
(2) For continuous variable, the search region in each round should be expanded from the center, which is the design variable for initial state.
(3) Accelerating technique for continuous variable incorporates decreasing the range of search region and the width of search tree as the update process proceeds. The width of search tree is based on the number of meshes generated in each round.
(4) Three kinds of accelerating techniques including geometric decay, linear decrease, and step reduction are also considered for continuous variable. Regardless of the range of search region and the width of search tree, geometric decay performs better than linear decrease and step reduction related to computational cost.
(5) For MVSMCTS formulation, various types of variables are considered in only one search tree.
(6) MDP framework, four steps of MCTS, UCB, the best reward, policy improvement, update process and accelerating technique for discrete variable and terminal condition in IMCTS formulation can also be employed in MVSMCTS formulation.

The proposed algorithm allows the agent to find an optimal solution. The numerical examples demonstrate that this algorithm provides results as good as metaheuristic algorithms in a reasonable time and applicable for practical engineering problems. In conclusion, this study suggests that MVSMCTS formulation is a powerful optimization technique for mixed variable structural optimization without tuning parameters.

The proposed method is expected to solve more complex problems such as simultaneous sizing, shape, and topology optimization of truss structures, which are our future research interests.